\documentstyle[epsf]{article}

\newcommand{\bigzerou}{%
\smash{\lower1.7ex\hbox{\bg 0}}}
\newtheorem{theorem}{Theorem}
\newtheorem{prop}{Proposition}
\newtheorem{defi}{Definition}
\newtheorem{cor}{Corollary}
\newtheorem{conj}{Conjecture}

\newtheorem{Rem}{Remark}

\newcommand{\ba}{\begin{eqnarray}}
\newcommand{\ea}{\end{eqnarray}}
\newcommand{\no}{\nonumber}

\def\d{{\partial}}
\newcommand{\mapright}[1]{%
\smash{\mathop{%
\hbox to 1.0cm{\rightarrowfill}}\limits^{#1}}}
\newcommand{\mapleft}[1]{%
\smash{\mathop{%
\hbox to 1.3cm{\leftarrowfill}}\limits^{#1}}}
\begin{document}
\title{
\begin{flushright}
  \begin{minipage}[b]{5em}
    \normalsize
    EPHOU 11-005
    ${}$      \\\\
  \end{minipage}
\end{flushright}
{\bf  Open Virtual Structure Constants and Mirror Computation of Open Gromov-Witten Invariants of Projective Hypersurfaces }}
\author{Masao Jinzenji ${}^{(1)}$, Masahide Shimizu ${}^{(2)}$ \\
\\
\it (1) Department of Mathematics, Graduate School of Science \\
\it Hokkaido University \\
\it  Kita-ku, Sapporo, 060-0810, Japan\\
{\it e-mail address: jin@math.sci.hokudai.ac.jp}\\
\it (2) Department of Physics, Graduate School of Science \\
\it Hokkaido University \\
\it  Kita-ku, Sapporo, 060-0810, Japan\\
{\it e-mail address: shimizu@particle.sci.hokudai.ac.jp}}
\maketitle
\begin{abstract}
In this paper, we generalize Walcher's computation of the open Gromov-Witten invariants of the quintic hypersurface to Fano and 
Calabi-Yau projective hypersurfaces. Our main tool is the open virtual structure constants. 
We also propose the generalized mirror transformation for the open Gromov-Witten invariants, 
some parts of which are proven explicitly. 
We also discuss possible modification of the multiple covering formula for the case of higher dimensional Calabi-Yau manifolds. 
The generalized disk invariants for some Calabi-Yau and Fano manifolds are shown and they are certainly integers after re-summation by the modified multiple covering formula. 
This paper also contains the direct integration method of the period integrals for higher dimensional Calabi-Yau hypersurfaces in the appendix. 
\end{abstract}
\section{Introduction}

Topological string theory and mirror symmetry are very powerful tools for attacking problems of enumerative geometry and Gromov-Witten theory. 
The first successful example is the celebrated work by Candelas et.al. \cite{CDGP}, 
where the number of rational curves of any degrees in the quintic $3$-fold was predicted. 
The result was rigorously proven later in various mathematical contexts, e.g., \cite{gi1,gi2,lly} by using localization computation \cite{kont}. 
Physically, the number of rational curves (and Gromov-Witten invariants) is related to the topological closed string amplitudes and 
extension to the open string sector was also initiated by physicists. 
In particular, 
the integral invariants related to the open string sector and the multiple covering formula were discussed in \cite{OV}. The 
disk invariants for several non-compact toric Calabi-Yau $3$-folds were predicted by using open mirror symmetry \cite{AV,AKV}, and 
also by using duality between topological string theory and  Chern-Simons theory \cite{LMV}. 
The disk invariants are number of holomorphic disks (thus, integer valued) and they are related to genus $0$, $1$-holed open Gromov-Witten invariants via the multiple covering formula. 
Earlier attempts to mathematical construction of the open Gromov-Witten invariant are contained in \cite{KL,L,GZ}. 
Recently, the study of open mirror symmetry for {\it compact} Calabi-Yau $3$-fold showed major breakthrough after appearance of the work by Walcher \cite{walcher}. 
In \cite{walcher}, the disk invariants of the quintic $3$-fold were computed by using both localization calculation based on \cite{solomon} and open mirror symmetry. 
These results were confirmed mathematically in \cite{psw}. 
After these works, there appeared many related works, e.g., the B-model side analyses \cite{mw,LLY2}, and 
the so-called off-shell idea \cite{JS,AHMM}. 
Prediction of the disk invariants for compact Calabi-Yau $3$-fold is now widely extended to various cases. 
For example, the disk invariants for several pfaffian Calabi-Yau varieties were predicted in \cite{SS} 
by using the method of direct integration of the period integrals developed in \cite{FNSS}. 
As a matter of course, mathematically rigorous construction of the open Gromov-Witten invariant have been tried in many works, 
e.g., \cite{ia,fukaya,PZ}. 

Meanwhile, 
some attempts to foundational understanding of mirror symmetry have been pursued by one of the present authors. 
In the series of his works, the notions of the virtual structure constant and the generalized mirror transformation were proposed. 
Roughly speaking, the virtual structure constants are the B-model analogue of the structure constants of small quantum cohomology ring. 
They were originally defined by certain recursive formulas \cite{CJ}. 
After \cite{CJ}, their relation to the B-model amplitudes and to the differential equation which governs the B-model geometry was clarified \cite{gm}. 
Later, it was observed that they can be evaluated by computing intersection numbers of the moduli space of polynomial maps with two marked points via localization computation \cite{vs,mmg}. 
The generalized mirror transformation translates these virtual structure constants into the ordinary structure constants, namely, the Gromov-Witten invariants. 
It was predicted in \cite{QCgeneral,gene} and proved mathematically in \cite{iri} as an
effect of coordinate change of the virtual Gauss-Manin system. 
It was also shown that these notions are very useful for not only Calabi-Yau manifolds but also Fano and general type manifolds \cite{QCgeneral,vGW}. 
In \cite{pmth}, the residue integral representations of Gromov-Witten invariants and the virtual structure constants are effectively  used so that 
we can reduce non-trivial relations coming from mirror symmetry predictions into simple algebraic identities of rational functions.

One of the main objectives of this paper is to extend these notions to the open string sector. 
Here, let us summarize the geometric settings used in this paper, 
which is mainly based on the treatment of the disk invariants in \cite{solomon,psw}. 
We consider Calabi-Yau (and Fano) manifolds 
which can be expressed as projective hypersurfaces (or complete intersections). 
In this paper, we consider the {\it real} defining equation, whose coefficients are all real numbers. More precisely, 
we consider the Fermat-type defining equation, 
\begin{equation}
\{ {X_1}^k+{X_2}^k+\dots +{X_N}^k=0 \}\subset CP^{N-1}, 
\label{defhy}
\end{equation}
where $k$ is an odd positive integer. We denote this hypersurface by $M_{N}^{k}$.
We also consider the special Lagrangian submanifold in it, 
which is defined as a fixed locus of the anti-holomorphic involution:
\begin{equation}
\{ {X_1}^k+{X_2}^k+\dots +{X_N}^k=0 \}\subset RP^{N-1}=:CP^{N-1}_{R}. 
\end{equation}
This locus turns out to be topologically $RP^{k-2}$ and have discrete $Z_2$ moduli. 
\begin{figure}[h]
      \epsfxsize=10cm
     \centerline{\epsfbox{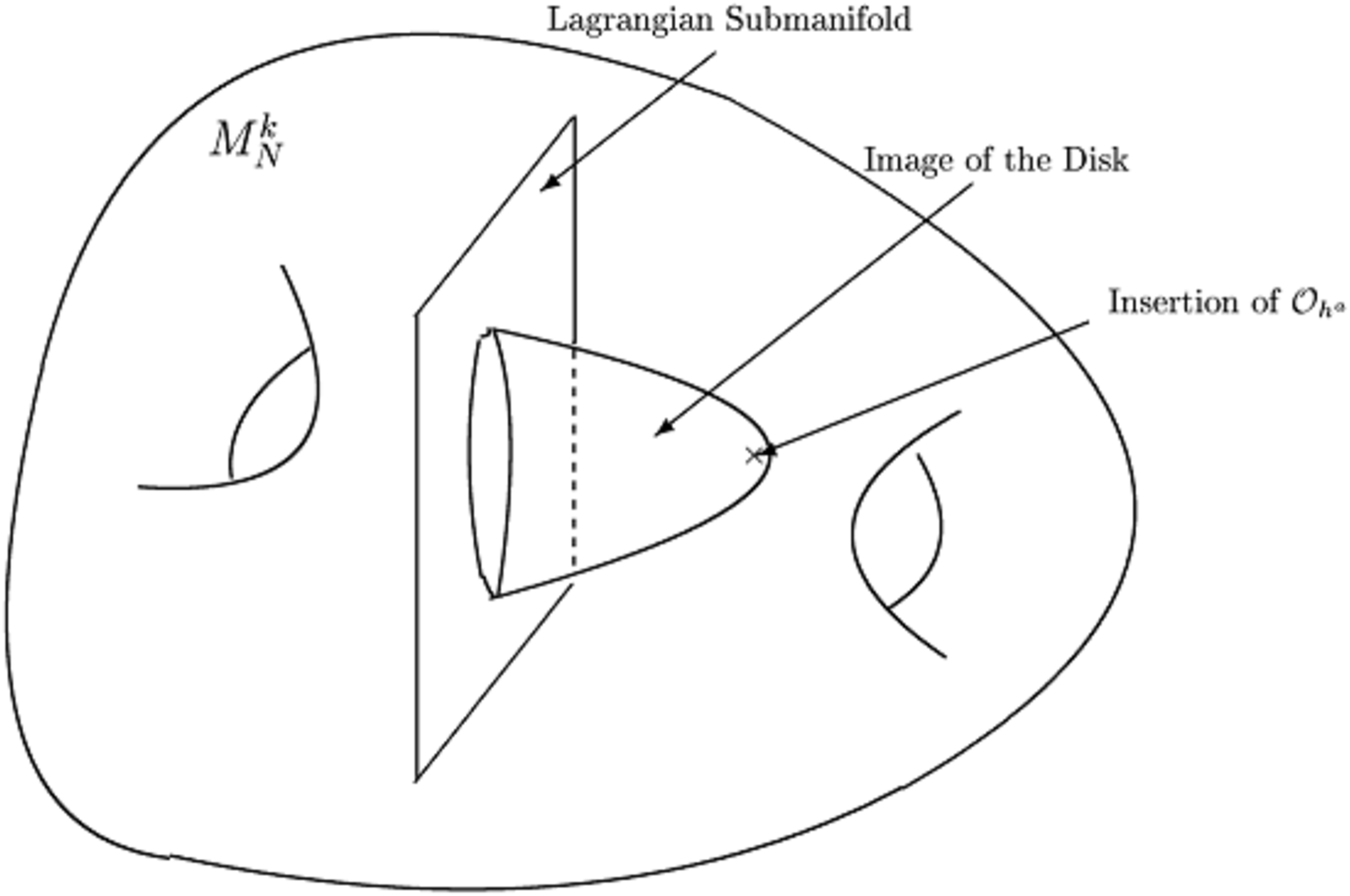}}
    \caption{\bf Geometrical Setting }
\label{open-sum}
\end{figure}
In this setting, we compute the 1-pointed open Gromov-Witten invariant $\langle{\cal O}_{h^{a}}\rangle_{disk,2d-1}$ 
where $h$ is the hyperplane class of $CP^{N-1}$. Of course, rigorous definition of $\langle{\cal O}_{h^{a}}\rangle_{disk,2d-1}$ 
is a hard problem from the point of view of mathematics. In this paper, we proceed by using the following heuristic argument. 
We expect that the number $\langle{\cal O}_{h^{a}}\rangle_{disk,2d-1}$ is obtained from integrating out the class $e(F_{2d-1})\wedge ev_{1}^{*}(h^a)$ on the moduli space of   
$\widetilde{M}_{D,1}(CP^{N-1}/CP^{N-1}_{R},2d-1)$ (we follow the notations in \cite{psw}). 
$\widetilde{M}_{D,1}(CP^{N-1}/CP^{N-1}_{R},2d-1)$ is the moduli space 
obtained from $\overline{M}_{D,1}(CP^{N-1}/CP^{N-1}_{R},2d-1)$, the moduli space of stable maps from 1-pointed  
disk to $CP^{N-1}$ (the boundary of a disk is mapped to $CP^{N-1}_{R}=RP^{N-1}$ and the one marked point is located 
inside the disk).
It is constructed  by generalizing the construction of 
$\overline{M}_{D,0}(CP^{4}/CP^{4}_{R},2d-1)$ in \cite{psw}. $F_{2d-1}$ is the vector bundle which guarantees that the image of the 
disk lies inside the hypersurface. By dimensional counting, we can see that $\langle{\cal O}_{h^{a}}\rangle_{disk,2d-1}$ is non-zero
only if $a=\frac{N-3+(N-k)(2d-1)}{2}$. Of course, $\frac{N-3+(N-k)(2d-1)}{2}$ must be an integer. 
From this heuristic definition, we can compute $\langle{\cal O}_{h^{a}}\rangle_{disk,2d-1}$ 
by generalizing the localization computation given in \cite{walcher}. We present some explicit formulas to compute $\langle{\cal O}_{h^{a}}\rangle_{disk,2d-1}$ for low degrees in Section \ref{PROOF}.
On the other hand, we define the open virtual structure constants, which are the B-model analogue of the open Gromov-Witten invariants. 
Our formula for the open virtual structure constant can be obtained by simple modification of the residue integral representation 
of the (closed) virtual structure constants. 
This modification comes from imitating the evaluation of open Gromov-Witten invariants via the localization computation in \cite{walcher}.
When $N=k$, i.e., the hypersurface in (\ref{defhy}) is Calabi-Yau, we can show that the generating function of the open virtual 
structure constants is obtained from the solutions of Picard-Fichs equation: 
\begin{eqnarray}
&&(\frac{d}{dx}-\frac{1}{2})\cdot\biggl((\frac{d}{dx})^{k-1}-k\cdot e^{x}\cdot
(k\frac{d}{dx}+k-1)\cdots
(k\frac{d}{dx}+2)\cdot(k\frac{d}{dx}+1)\biggr)w(x)=0.
\label{openhyper} 
\end{eqnarray}  
We also propose the generalized mirror transformation that relates the open virtual structure constants to $\langle{\cal O}_{h^{a}}\rangle_{disk,2d-1}$. The proof of the generalized mirror transformation for low degrees is given in Section \ref{PROOF}. 
  
With these results, we perform the B-model computation of $\langle{\cal O}_{h^{\frac{k-3}{2}}}\rangle_{disk,2d-1}$ for 
the Calabi-Yau hypersurface $M_{k}^{k}$. Then we propose the multiple covering formula for $\langle{\cal O}_{h^{\frac{k-3}{2}}}\rangle_{disk,2d-1}$ and present the numerical results of the corresponding integral disk invariants up to 
$k=13$. We also present some numerical results for Calabi-Yau complete intersections and Fano hypersurfaces. 
Furthermore, we discuss the direct integration method of the period integrals and the chain integrals, which is a natural extension of the method discussed in \cite{FNSS,SS} to higher dimensional cases. 
We think that 
our B-model computation of the disk invariants of higher dimensional Calabi-Yau manifold is new both in physics and mathematics.     

Organization of this paper is the following. 
First, in Section \ref{B}, we focus on the B-model side and propose the closed formula of the generating function of $\langle{\cal O}_{h^{\frac{k-3}{2}}}\rangle_{disk,2d-1}$ for Calabi-Yau hypersurfaces, which is written in terms of 
the ordinary virtual structure constants and the solution of the inhomogeneous Picard-Fuchs equation. 
In Section \ref{OVSC}, the definition of the open virtual structure constants are presented and its relation to the arguments in Section \ref{B} is discussed. 
We also propose the generalized mirror transformation for the open Gromov-Witten invariants of $M_{N}^{k}$. 
In Section \ref{DATA}, the disk invariants for various Calabi-Yau and Fano hypersurfaces/complete intersections are exhibited. 
In Calabi-Yau cases, these invariants are obtained by adopting the formulas developed in the previous sections and the modified multiple covering formula. 
In Section \ref{PROOF}, we introduce the residue integral representation of the open Gromov-Witten invariants and prove the generalized mirror transformation proposed in Section \ref{OVSC} up to degree $5$ explicitly. 
Appendix \ref{DIRECT} is devoted to discussions of the direct integration method of the period integrals and the chain integrals, 
by which we can compute the fundamental period and the domainwall tension for Calabi-Yau projective hypersurfaces.

{\bf Acknowledgment} M.J. would like to thank Miruko Jinzenji for kind encouragement. M.S. would like to thank Dr. Yutaka Tobita for 
assistance of numerical computations using Mathematica.  Research of M.J. is partially supported by JSPS grant No. 22540061. 
\section{B-model Computation}\label{B}
In this section, we present a conjecture that generalizes Walcher's computation of the disk invariants (genus $0$, $1$-holed open Gromov-Witten invariants) of the quintic 
hypersurface in $CP^{4}$ to the degree $k$ Calabi-Yau hypersurface in $CP^{k-1}$ ($M_{k}^{k}$).  
In this paper, we assume that the positive integer $k$ always takes odd value.
First, we introduce the virtual structure constants $\tilde{L}_{n}^{k,k,d}$ ($n=0$, $1$, $\cdots$, $k-1$) 
for $M_{k}^{k}$, which was introduced in \cite{gm}.
It is determined by the following equation:  
\begin{eqnarray}
&&\biggl((\frac{d}{dx})^{k-1}-k\cdot e^{x}\cdot
(k\frac{d}{dx}+k-1)\cdots
(k\frac{d}{dx}+2)\cdot(k\frac{d}{dx}+1)\biggr)w(x)\no\\
&&=\frac{1}{\tilde{L}^{k,k}_{k-1}(e^x)}(\frac{d}{dx}\frac{1}{\tilde{L}^{k,k}_{k-2}(e^x)}(\frac{d}{dx}
\cdots \frac{1}{\tilde{L}^{k,k}_{1}(e^x)}(\frac{d}{dx}\frac{w(x)}{\tilde{L}^{k,k}_{0}(e^x)})\cdots)),
\label{hyper} 
\end{eqnarray}
where $\tilde{L}_{n}^{k,k}(e^x)$ is the generating function of $\tilde{L}_{n}^{k,k,d}$:
\begin{equation}
\tilde{L}_{n}^{k,k}(e^x):=1+\sum_{d=1}^{\infty}\tilde{L}_n^{k,k,d}e^{dx}.
\label{vir1}
\end{equation}
In (\ref{hyper}), $w(x)$ is an arbitrary function with adequate differentiable property. 
We can construct $\tilde{L}_{n}^{k,k}(e^x)$ that satisfies (\ref{hyper}) from the solutions of the differential 
equation:
\begin{equation}
\biggl((\frac{d}{dx})^{k-1}-k\cdot e^{x}\cdot
(k\frac{d}{dx}+k-1)\cdots
(k\frac{d}{dx}+2)\cdot(k\frac{d}{dx}+1)\biggr)w(x)=0.
\label{hyper2}
\end{equation}
Linearly independent solutions of (\ref{hyper2}) around $x=-\infty$ are explicitly given by $w_{j}(x)$ ($j=0$, $1$, $2$, $\cdots$, $k-2$): 
\begin{eqnarray}
w(x,y)&:=&\sum_{d=0}^{\infty}\frac{\prod_{j=1}^{kd}(j+ky)}{\prod_{j=1}^{d}(j+y)^k}e^{(d+y)x}.\no\\
w_j(x)&:=&\frac{1}{j!}\frac{\d^{j}w}{\d y^{j}}(x,0).
\label{sol}
\end{eqnarray}
Then $\tilde{L}_{n}^{k,k}(e^x)$ is inductively determined by the following relation\footnote{
In (\ref{recur}), we need to use formally $w_{k-1}(x)$ to determine $\tilde{L}_{k-1}^{k,k}(e^x)$ though it is 
not a solution of (\ref{hyper2}). }: 
\begin{eqnarray}
\tilde{L}_{0}^{k,k}(x)&=&w_{0}(x),\no\\
\tilde{L}_{j}^{k,k}(e^x)&=&\frac{d}{dx}(\frac{1}{\tilde{L}^{k,k}_{j-1}(e^x)}\frac{d}{dx}(\frac{1}{\tilde{L}^{k,k}_{j-2}(e^x)}
\frac{d}{dx}(\frac{1}{\tilde{L}^{k,k}_{j-3}(e^x)}\cdots\frac{d}{dx}(\frac{1}{\tilde{L}^{k,k}_{1}(e^x)}\frac{d}{dx}\frac{w_j(x)}{\tilde{L}_{0}^{k,k}(e^x)})
\cdots))).
\label{recur}
\end{eqnarray}
For later use, we also introduce the residue integral representation \cite{mmg} of the virtual structure constant $\tilde{L}_{n}^{k,k,d}$ as follows:
\begin{eqnarray}
\tilde{L}_n^{k,k,d}:=&&\frac{d}{k}\cdot\frac{1}{(2\pi\sqrt{-1})^{d+1}}\oint_{C_{d}}\frac{dz_{d}}{(z_{d})^k}
\oint_{C_{d-1}}\frac{dz_{d-1}}{(z_{d-1})^k}
\cdots\oint_{C_{0}}\frac{dz_{0}}{(z_{0})^k}(z_{0})^{k-2-n}\bigl(\prod_{i=1}^{d}e^k(z_{i-1},z_i)\bigr)\times\no\\
&&\times\bigl(\prod_{i=1}^{d-1}\frac{1}{kz_i(2z_i-z_{i-1}-z_{i+1})}\bigr)(z_d)^{n-1},
\label{def2}
\end{eqnarray}
where $\frac{1}{2\pi\sqrt{-1}}\oint_{C_{i}}$ means taking residues at $z_i=0$, $z_{i}=\frac{z_{i-1}+z_{i+1}}{2}$ if $i=1$, $\cdots$, $d-1$ and at $z_i=0$
if $i=0$, $d$. $e^{k}(z,w)$ is a polynomial in $z$ and $w$, which is given by,
\begin{eqnarray}
e^{k}(z,w):=\prod_{j=0}^{k}(jz+(k-j)w). \label{e^k}
\end{eqnarray}
Note that, instead of (\ref{e^k}), it is possible to use the following rather complicated rational function: 
\begin{equation}
{E}^{N,k}_d(z,w):=\frac{\prod_{j=0}^{k d} (\frac{j z+(k d-j)w}{d})}{\prod_{j=1}^{d-1} (\frac{j w+(d-j) z}{d})^N}\label{E^k}. 
\label{ENK}
\end{equation}
Using (\ref{ENK}), we can write down an alternate residue integral representation of $\tilde{L}_n^{k,k,d}$:
\begin{eqnarray}
\tilde{L}_{n}^{k,k,d}&=&\frac{d}{k}
\sum_{\sigma_{d}\in OP_{d}}\frac{1}{(2\pi\sqrt{-1})^{l(\sigma_{d})+1}}\oint_{C_{(0)}} \frac{dz_{0}}{(z_{0})^k}\cdots
\oint_{C_{(0)}} \frac{d z_{l(\sigma_{d})} }{(z_{l(\sigma_{d})})^k}(z_{0})^{k-2-n}(z_{l(\sigma_{d})})^{n-1}\times\no\\
&&\times\prod_{j=1}^{l(\sigma_{d})-1}\frac{1}
{\biggl( \frac{z_{j}-z_{j-1}}{d_{j}}+\frac{z_{j}-z_{j+1}}{d_{j+1}}\biggr)kz_{j}}\prod_{j=1}^{l(\sigma_{d})}
\frac{E^{k,k}_{d_j}(z_{j-1},z_j)}{d_j},
\label{int}
\end{eqnarray}
where $OP_{d}$ is the set of ordered partitions of the positive integer $d$:
\begin{eqnarray}
OP_d:=\{\sigma_{d}:=(d_1,d_2,\cdots,d_{l(\sigma_d)})\;|\;d_i\geq 0,\;\sum_{i=1}^{l(\sigma_d)}d_i=d\;\}.
\label{op}
\end{eqnarray}
$\frac{1}{2\pi\sqrt{-1}}\oint_{C_{(0)}}dz$ means operation of taking residue at $z=0$.
Since we can derive (\ref{int}) from (\ref{def2}) by taking residues at $z_{i}=\frac{z_{i-1}+z_{i+1}}{2}$ first,    
the resulting $\tilde{L}_n^{k,k,d}$ is completely the same. See \cite{mmg} for detailed derivation.
In this paper, we choose the simpler formula (\ref{def2}). 
The formula (\ref{int}) is needed when we consider the open string modification of the virtual structure constants. 

Next, we introduce the function:
\begin{equation}
\tau_{k}(x):=\sum_{d=1}^{\infty}2\frac{(k(2d-1))!!}{((2d-1)!!)^{k}}e^{\frac{2d-1}{2}x},
\label{superpotential}
\end{equation}
that satisfies the inhomogeneous version of (\ref{hyper2}) as follows: 
\begin{equation}
\biggl((\frac{d}{dx})^{k-1}-k\cdot e^{x}\cdot
(k\frac{d}{dx}+k-1)\cdots
(k\frac{d}{dx}+2)\cdot(k\frac{d}{dx}+1)\biggr)\tau_{k}(x)=C\cdot\exp(\frac{1}{2}x),
\label{hyper3}
\end{equation}
where $C$ is some constant. 
This kind of inhomogeneity was first considered in \cite{walcher} and precise value of the constant for the quintic $3$-fold case was determined in \cite{mw}. 
Further discussions are contained in \cite{LLY2}. 
An alternative method to obtain (\ref{superpotential}) was proposed in \cite{FNSS,SS}. 
In Appendix \ref{DIRECT}, we discuss this method and apply it to the higher dimensional cases. 
With this set-up, we introduce the following function: 
\begin{eqnarray}
&&F_{o}^{k}(x):=2^{\frac{k-3}{2}}\frac{1}{\tilde{L}^{k,k}_{\frac{k-3}{2}}(e^{x})}\frac{d}{dx}\frac{1}{\tilde{L}^{k,k}_{\frac{k-5}{2}}(e^{x})} 
\cdots\frac{d}{dx}\frac{1}{\tilde{L}^{k,k}_{1}(e^{x})} 
\frac{d}{dx}\frac{\tau_{k}(x)}{\tilde{L}^{k,k}_{0}(e^{x})} .\no\\
\end{eqnarray}
\begin{conj}
By using the mirror map used in the mirror computation of $M_{k}^{k}$:
\begin{equation}
t(x):=x+\sum_{d=1}^{\infty}\frac{\tilde{L}_{1}^{k,k,d}}{d}e^{dx},\label{mirror-map}
\end{equation}
$F_{o}^{k}(x)$ gives us the generating function of the open Gromov-Witten invariants of $M_{k}^{k}$ via the following formula:
\begin{equation}
F_{o}^{k}(x(t))=\sum_{d=1}^{\infty}\langle{\cal O}_{h^{\frac{k-3}{2}}}\rangle_{disk,2d-1}e^{\frac{2d-1}{2}t}. 
\end{equation}
Here, $x(t)$ is the inverse map of (\ref{mirror-map}). 
\end{conj}
\section{Open Virtual Structure Constants}\label{OVSC}
In this section, we consider the open Gromov-Witten invariants for $M_{N}^{k}$. 
First, we briefly explain how we reached Conjecture 1. 
It is based on the A-model computation of the open Gromov-Witten invariants. 
We introduce here the following rational function in $z$ for  positive integer $d$:
\begin{eqnarray}
f^{N,k}_{2d-1}(z)&:=&\frac{2}{2d-1}\cdot\frac{\displaystyle{\prod_{j=0}^{kd-\frac{k+1}{2}}(\frac{j(-z)+(k(2d-1)-j)z}{2d-1})}}
{\displaystyle{\prod_{j=1}^{d-1}(\frac{j(-z)+(2d-1-j)z}{2d-1})^N}}.
\label{def1}
\end{eqnarray}
By using (\ref{def1}), we define the open virtual structure constant as follows. 
\begin{defi} For positive integer $d$, we define,
\begin{eqnarray}
w_{disk}^{N,k}({\cal O}_{h^a})_{2d-1}&&:=
\sum_{j=0}^{d-1}\frac{1}{(2\pi\sqrt{-1})^{j+1}}\oint_{C_{j}}\frac{dz_{j}}{(z_{j})^N}
\oint_{C_{j-1}}\frac{dz_{j-1}}{(z_{j-1})^N}
\cdots\oint_{C_{0}}\frac{dz_{0}}{(z_{0})^N}f_{2d-2j-1}(z_0)\times\no\\
&&\times\bigl(\prod_{i=1}^{j}e^k(z_{i-1},z_i)\bigr)
\frac{1}{kz_0(\frac{2}{2d-2j-1}z_0+z_0-z_{1})}\bigl(\prod_{i=1}^{j-1}\frac{1}{kz_i(2z_i-z_{i-1}-z_{i+1})}\bigr)(z_j)^a,
\label{def3}
\end{eqnarray}
where $\oint_{C_{i}}$ means taking residues at $z_i=0,\frac{z_{i-1}+z_{i+1}}{2}$ if $i=1$, $\cdots$, $j-1$ and at $z_i=0$
if $i=0$, $j$. We integrate the variable $z_{i}$ in ascending order of the subscript $i$. 
\end{defi}
Let us explain the geometrical meaning of the formula (\ref{def3}). In \cite{walcher}, 
$0$-point open Gromov-Witten invariant of degree $2d-1$ of the quintic hypersurface in $CP^4$ 
is evaluated by integrating out the square root of the class $c_{top}(R^{0}\pi_{*}ev_{1}^{*}({\cal O}_{CP^4}(5)))$
on the moduli space $\widetilde{M}_{D,0}(CP^4/CP^4_{R},2d-1)$ (we follow the notation in \cite{psw}). 
In our case, instead of considering $0$-point invariants, 
we evaluate $1$-point invariants $\langle{\cal O}_{h^a}\rangle_{disk,2d-1}$ 
by integrating out the 
class $e(F_d)\wedge ev_{1}^{*}(h^a)$ on the moduli space of   
$\widetilde{M}_{D,1}(CP^{N-1}/CP^{N-1}_{R},2d-1)$. $e(F_d)$ is the generalization of the class 
$\sqrt{c_{top}(R^{0}\pi_{*}ev_{1}^{*}({\cal O}_{CP^4}(5)))}$, which guarantees that the image of the (stable) disk 
lies inside $M_{N}^{k}$.
To understand (\ref{def3}), let us remember the situation for the closed string case. 
In \cite{mmg}, we introduced the intersection number 
$w({\cal O}_{h^a}{\cal O}_{h^b})_{0,d}$ of the moduli space of polynomial maps from $CP^1$ to $CP^{N-1}$ 
with two marked points (we denote it by $\widetilde{Mp}_{0,2}(CP^{N-1},d)$). $w({\cal O}_{h^a}{\cal O}_{h^b})_{0,d}$ 
corresponds to the B-model analogue of the Gromov-Witten invariant $\langle{\cal O}_{h^a}{\cal O}_{h^b}\rangle_{0,d}$.  
It is evaluated by integrating out the class $c_{top}({\cal E}_{d}^{k})\wedge ev_{1}^{*}(h^a)\wedge ev_{2}^{*}(h^b)$
on $\widetilde{Mp}_{0,2}(CP^{N-1},d)$
that corresponds to $c_{top}(R^{0}\pi_{*}ev_{3}^{*}({\cal O}_{CP^{N-1}}(k))) \wedge ev_{1}^{*}(h^a)\wedge ev_{2}^{*}(h^b)$
on $\overline{M}_{0,2}(CP^{N-1},d)$, the moduli space of stable maps.  
By applying the localization theorem, it is explicitly given by the following formula.
\begin{eqnarray}
w({\cal O}_{h^a}{\cal O}_{h^b})_{0,d} &=&
\sum_{\sigma_{d}\in OP_{d}}\frac{1}{(2\pi\sqrt{-1})^{l(\sigma_{d})+1}}\oint_{C_{(0)}} \frac{dz_{0}}{(z_{0})^N}\cdots
\oint_{C_{(0)}} \frac{d z_{l(\sigma_{d})} }{(z_{l(\sigma_{d})})^N}(z_{0})^{N-2-n}(z_{l(\sigma_{d})})^{n-1+(N-k)d}\times\no\\
&&\times\prod_{j=1}^{l(\sigma_{d})-1}\frac{1}
{\biggl( \frac{z_{j}-z_{j-1}}{d_{j}}+\frac{z_{j}-z_{j+1}}{d_{j+1}}\biggr)kz_{j}}\prod_{j=1}^{l(\sigma_{d})}
\frac{E^{k,k}_{d_j}(z_{j-1},z_j)}{d_j}\no\\
\bigl(&=&\frac{1}{(2\pi\sqrt{-1})^{d+1}}\oint_{C_{d}}\frac{dz_{d}}{(z_{d})^N}
\oint_{C_{d-1}}\frac{dz_{d-1}}{(z_{d-1})^N}
\cdots\oint_{C_{0}}\frac{dz_{0}}{(z_{0})^N}(z_{0})^a\bigl(\prod_{i=1}^{d}e^k(z_{i-1},z_i)\bigr)\times\no\\
&&\times\bigl(\prod_{i=1}^{d-1}\frac{1}{kz_i(2z_i-z_{i-1}-z_{i+1})}\bigr)(z_d)^{b}\bigr).
\label{wint}
\end{eqnarray}
In the last line of (\ref{wint}), $\frac{1}{2\pi\sqrt{-1}}\oint_{C_i}$ has the same meaning as the one in (\ref{def2}).
With this formula in mind, we consider the open string analogue of $w({\cal O}_{h^a}{\cal O}_{h^b})_{0,d}$, say, 
$w_{disk}^{N,k}({\cal O}_{h^a})_{2d-1}$, by integrating out $\sqrt{c_{top}({\cal E}_d^{k})}ev_{1}^{*}(h^a)$ 
on $\widetilde{Mp}_{D,1}(CP^{N-1}/CP^{N-1}_{R},2d-1)$ , the moduli space of polynomial maps with $Z_{2}$ symmetry 
compatible with the anti-holomorphic involution of $CP^{N-1}$. 
By modifying the localization method used in \cite{vs} according to the discussion given in \cite{psw}, we propose 
the following definition of $w_{disk}^{N,k}({\cal O}_{h^a})_{2d-1}$: 
\begin{eqnarray}
w_{disk}^{N,k}({\cal O}_{h^a})_{2d-1}&&:=
\sum_{j=0}^{d-1}\sum_{\sigma_j\in OP_j}\frac{1}{(2\pi\sqrt{-1})^{l(\sigma_j)+1}}\oint_{C_{(0)}}\frac{dz_{l(\sigma_j)}}{(z_{l(\sigma_j)})^N}
\oint_{C_{(0)}}\frac{dz_{l(\sigma_j)-1}}{(z_{l(\sigma_j)-1})^N}
\cdots\oint_{C_{(0)}}\frac{dz_{0}}{(z_{0})^N}f_{2d-2j-1}(z_0)\times\no\\
&&\times\bigl(\prod_{i=1}^{l(\sigma_j)}\frac{E^{N,k}_{j_i}(z_{i-1},z_i)}{j_i}\bigr)
\frac{1}{kz_0(\frac{2}{2d-2j-1}z_0+\frac{z_0-z_{1}}{j_1})}\bigl(\prod_{i=1}^{l(\sigma_j)-1}\frac{1}{kz_i(\frac{z_i-z_{i-1}}{j_i}+\frac{z_i-z_{i+1}}{j_{i+1}})}\bigr)(z_{l(\sigma_j)})^a. 
\label{odef3}
\end{eqnarray}
$f_{2d-2j-1}(z_0)$ in (\ref{odef3}) can be interpreted as the modified version of $E^{N,k}_{2d-2j-1}(z,w)$ in (\ref{E^k}), where modification is just to set $-w=z=z_0$, to take square root and to multiply the factor $\frac{2}{2d-2j-1}$. 
The factor $\frac{2}{2d-2j-1}$ comes from 
the action of the discrete automorphism group that acts on the edge with degree $\frac{2d-2j-1}{2}$.
The factor $\displaystyle{\frac{1}{(\frac{2}{2d-2j-1}z_0+\frac{z_0-z_{1}}{j_1})}}$ in (\ref{odef3}) comes from the equality:
\begin{equation}
\frac{1}{(\frac{2}{2d-2j-1}z_0+\frac{z_0-z_{1}}{j_1})}=\frac{1}{(\frac{z_0-(-z_0)}{2d-2j-1}+\frac{z_0-z_{1}}{j_1})}.
\end{equation}  
We can reduce the formula (\ref{odef3}) to (\ref{def3}) by using the same trick which reduces the formula (\ref{int}) to (\ref{def2}).
The resulting formula (\ref{def3}) can be interpreted as the sum of the contributions associated with the graph in Figure 2\footnote{In the figures in this paper, the degrees associated with the edges of the graphs are half of the degrees used in the body of this paper.}.
\begin{figure}[h]
      \epsfxsize=12cm
     \centerline{\epsfbox{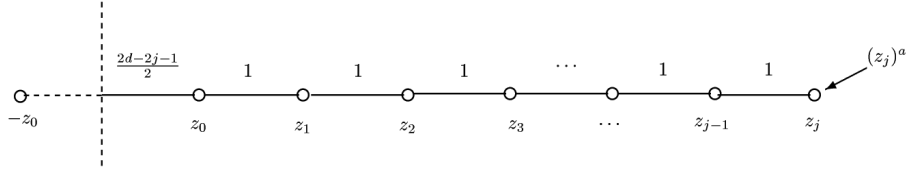}}
    \caption{\bf The graph that contributes to $w_{disk}^{N,k}({\cal O}_{h^a})_{2d-1}$}
\label{g4}
\end{figure}

Now, let us  present a key theorem that leads us to Conjecture 1.
\begin{theorem} Let us assume that $N=k$. If we expand $F_{o}^{k}(x)$ in the form:
\begin{equation}
F_{o}^{k}(x)=\sum_{d=1}^{\infty}c_{2d-1}^{k}e^{(d-\frac{1}{2})x},
\end{equation}
the following equality holds.
\begin{equation}
c_{2d-1}^{k}=w_{disk}^{k,k}({\cal O}_{h^{\frac{k-3}{2}}})_{2d-1}.\\
\end{equation}
\end{theorem}
{\it proof)}
First, we introduce the following function. 
\begin{equation}
G_{o}^{k}(x):=\sum_{d=1}^{\infty}w_{disk}^{k,k}({\cal O}_{h^{\frac{k-3}{2}}})_{2d-1}e^{(d-\frac{1}{2})x}.
\end{equation}
Then it is enough for us to show the following equality.
\begin{eqnarray}
\tau_{k}(x)=\frac{1}{2^{\frac{k-3}{2}}}\tilde{L}_{0}^{k,k}(e^x)\int_{-\infty}^{x}dx_1 \tilde{L}_{1}^{k,k}(e^{x_1})\int_{-\infty}^{x_1}dx_2 \tilde{L}_{2}^{k,k}(e^{x_{2}})\cdots
\int_{-\infty}^{x_{\frac{k-5}{2}}}dx_{\frac{k-3}{2}} \tilde{L}_{\frac{k-3}{2}}^{k,k}(e^{x_{\frac{k-3}{2}}})G_{o}^{k}(x_{\frac{k-3}{2}}).
\end{eqnarray}
First, we compute $\frac{1}{2}\int_{-\infty}^{x_{\frac{k-5}{2}}}dx_{\frac{k-3}{2}} \tilde{L}_{\frac{k-3}{2}}^{k,k}(e^{x_{\frac{k-3}{2}}})G_{o}^{k}(x_{\frac{k-3}{2}})$. 
The coefficient of $e^{(d-\frac{1}{2})x}$ in $\tilde{L}_{\frac{k-3}{2}}^{k,k}(e^{x})G_{o}^{k}(x)$
is given by,
\begin{eqnarray}
&&\sum_{j=0}^{d-1}w_{disk}^{k,k}({\cal O}_{h^{\frac{k-3}{2}}})_{2d-2j-1}\tilde{L}_{\frac{k-3}{2}}^{k,k,j}\no\\
&&=\sum_{j=0}^{d-1}
\Big[\sum_{l=0}^{d-j-1}\frac{1}{(2\pi\sqrt{-1})^{l+1}}\oint_{C_{l}}\frac{dz_{l}}{(z_{l})^k}
\oint_{C_{l-1}}\frac{dz_{l-1}}{(z_{l-1})^k}
\cdots\oint_{C_{0}}\frac{dz_{0}}{(z_{0})^k}f_{2d-2j-2l-1}(z_0)\bigl(\prod_{i=1}^{l}e^k(z_{i-1},z_i)\bigr)\times\no\\
&&\times\frac{1}{kz_0(\frac{2}{2d-2j-2l-1}z_0+z_0-z_{1})}\bigl(\prod_{i=1}^{l-1}\frac{1}{kz_i(2z_i-z_{i-1}-z_{i+1})}\bigr)(z_l)^{\frac{k-3}{2}}\Big]\times
\no\\
&&\times\Big[\frac{j}{k}\cdot\frac{1}{(2\pi\sqrt{-1})^{j+1}}\oint_{C_{j}}\frac{dz_{j}}{(z_{j})^k}
\oint_{C_{j-1}}\frac{dz_{j-1}}{(z_{j-1})^k}
\cdots\oint_{C_{0}}\frac{dz_{0}}{(z_{0})^k}(z_{0})^{k-2-\frac{k-3}{2}}\bigl(\prod_{i=1}^{j}e^k(z_{i-1},z_i)\bigr)\times\no\\
&&\times\bigl(\prod_{i=1}^{j-1}\frac{1}{kz_i(2z_i-z_{i-1}-z_{i+1})}\bigr)(z_j)^{\frac{k-3}{2}-1}\Big].
\end{eqnarray}
At this stage, we can rewrite the product of two residue integrals into one residue integral 
if $l\neq 0$, $j\neq 0$. 
\begin{eqnarray}
&&\frac{1}{(2\pi\sqrt{-1})^{l+1}}\oint_{C_{l}}\frac{dz_{l}}{(z_{l})^k}
\oint_{C_{l-1}}\frac{dz_{l-1}}{(z_{l-1})^k}
\cdots\oint_{C_{0}}\frac{dz_{0}}{(z_{0})^k}f_{2d-2j-2l-1}(z_0)\bigl(\prod_{i=1}^{l}e^k(z_{i-1},z_i)\bigr)\times\no\\
&&\times\frac{1}{kz_0(\frac{2}{2d-2j-2l-1}z_0+z_0-z_{1})}\bigl(\prod_{i=1}^{l-1}\frac{1}{kz_i(2z_i-z_{i-1}-z_{i+1})}\bigr)(z_l)^{\frac{k-3}{2}}\times
\no\\
&&\times\frac{j}{k}\cdot\frac{1}{(2\pi\sqrt{-1})^{j+1}}\oint_{C_{j}}\frac{dz_{j}}{(z_{j})^k}
\oint_{C_{j-1}}\frac{dz_{j-1}}{(z_{j-1})^k}
\cdots\oint_{C_{0}}\frac{dz_{0}}{(z_{0})^k}(z_{0})^{k-2-\frac{k-3}{2}}\bigl(\prod_{i=1}^{j}e^k(z_{i-1},z_i)\bigr)\times\no\\
&&\times\bigl(\prod_{i=1}^{j-1}\frac{1}{kz_i(2z_i-z_{i-1}-z_{i+1})}\bigr)(z_j)^{\frac{k-3}{2}-1}\no\\
&=&\frac{1}{(2\pi\sqrt{-1})^{l+j+1}}\oint_{C_{l+j}}\frac{dz_{l+j}}{(z_{l+j})^k}
\oint_{C_{l+j-1}}\frac{dz_{l+j-1}}{(z_{l+j-1})^k}
\cdots\oint_{C_{0}}\frac{dz_{0}}{(z_{0})^k}f_{2d-2j-2l-1}(z_0)\bigl(\prod_{i=1}^{l}e^k(z_{i-1},z_i)\bigr)\times\no\\
&&\times\frac{1}{kz_0(\frac{2}{2d-2j-2l-1}z_0+z_0-z_{1})}\bigl(\prod_{i=1}^{l-1}\frac{1}{kz_i(2z_i-z_{i-1}-z_{i+1})}\bigr)\times
\no\\
&&\times\frac{j}{kz_l}\cdot\bigl(\prod_{i=l+1}^{j+l}e^k(z_{i-1},z_i)\bigr)
\bigl(\prod_{i=j+1}^{j+l-1}\frac{1}{kz_i(2z_i-z_{i-1}-z_{i+1})}\bigr)(z_{j+l})^{\frac{k-3}{2}-1}\no\\
&=&\frac{1}{(2\pi\sqrt{-1})^{l+j+1}}\oint_{C_{l+j}}\frac{dz_{l+j}}{(z_{l+j})^k}
\oint_{C_{l+j-1}}\frac{dz_{l+j-1}}{(z_{l+j-1})^k}
\cdots\oint_{C_{0}}\frac{dz_{0}}{(z_{0})^k}f_{2d-2j-2l-1}(z_0)\bigl(\prod_{i=1}^{j+l}e^k(z_{i-1},z_i)\bigr)\times\no\\
&&\times\frac{1}{kz_0(\frac{2}{2d-2j-2l-1}z_0+z_0-z_{1})}\bigl(\prod_{i=1}^{l+j-1}\frac{1}{kz_i(2z_i-z_{i-1}-z_{i+1})}\bigr)
(j(2z_{l}-z_{l-1}-z_{l+1}))(z_{j+l})^{\frac{k-3}{2}-1}.
\end{eqnarray}
This equality leads us to,  
\begin{eqnarray}
&&\sum_{j=0}^{d-1}w_{disk}^{k,k}({\cal O}_{h^{\frac{k-3}{2}}})_{2d-2j-1}\tilde{L}_{\frac{k-3}{2}}^{k,k,j}\no\\
&=&\sum_{m=0}^{d-1}\frac{1}{(2\pi\sqrt{-1})^{m+1}}\oint_{C_{m}}\frac{dz_{m}}{(z_{m})^k}
\oint_{C_{m-1}}\frac{dz_{m-1}}{(z_{m-1})^k}
\cdots\oint_{C_{0}}\frac{dz_{0}}{(z_{0})^k}f_{2d-2m-1}(z_0)\bigl(\prod_{i=1}^{m}e^k(z_{i-1},z_i)\bigr)\times\no\\
&&\times\frac{1}{kz_0(\frac{2}{2(d-m)-1}z_0+z_0-z_{1})}\bigl(\prod_{i=1}^{m-1}\frac{1}{kz_i(2z_i-z_{i-1}-z_{i+1})}\bigr)\times\no\\
&&\times(m(\frac{2}{2(d-m)-1}z_0+z_0-z_{1})+\sum_{j=1}^{m-1}j(2z_{m-j}-z_{m-j-1}-z_{m-j+1})+z_m)(z_{m})^{\frac{k-3}{2}-1}\no\\
&=&\sum_{m=0}^{d-1}\frac{1}{(2\pi\sqrt{-1})^{m+1}}\oint_{C_{m}}\frac{dz_{m}}{(z_{m})^k}
\oint_{C_{m-1}}\frac{dz_{m-1}}{(z_{m-1})^k}
\cdots\oint_{C_{0}}\frac{dz_{0}}{(z_{0})^k}f_{2d-2m-1}(z_0)\bigl(\prod_{i=1}^{m}e^k(z_{i-1},z_i)\bigr)\times\no\\
&&\times\frac{1}{kz_0(\frac{2}{2(d-m)-1}z_0+z_0-z_{1})}\bigl(\prod_{i=1}^{m-1}\frac{1}{kz_i(2z_i-z_{i-1}-z_{i+1})}\bigr)
\frac{d-\frac{1}{2}}{d-m-\frac{1}{2}}z_0(z_{m})^{\frac{k-3}{2}-1}.
\end{eqnarray} 
Therefore, we have the following equality.
\begin{eqnarray}
&&\frac{1}{2}\int_{-\infty}^{x_{\frac{k-5}{2}}}dx_{\frac{k-3}{2}} \tilde{L}_{\frac{k-3}{2}}^{k,k}(e^{x_{\frac{k-3}{2}}})G_{o}^{k}(x_{\frac{k-3}{2}})\no\\
&=&\sum_{d=1}^{\infty}e^{(d-\frac{1}{2})x_{\frac{k-5}{2}}}
\sum_{m=0}^{d-1}\frac{1}{(2\pi\sqrt{-1})^{m+1}}\oint_{C_{m}}\frac{dz_{m}}{(z_{m})^k}
\oint_{C_{m-1}}\frac{dz_{m-1}}{(z_{m-1})^k}
\cdots\oint_{C_{0}}\frac{dz_{0}}{(z_{0})^k}f_{2d-2m-1}(z_0)\bigl(\prod_{i=1}^{m}e^k(z_{i-1},z_i)\bigr)\times\no\\
&&\times\frac{1}{kz_0(\frac{2}{2(d-m)-1}z_0+z_0-z_{1})}\bigl(\prod_{i=1}^{m-1}\frac{1}{kz_i(2z_i-z_{i-1}-z_{i+1})}\bigr)
\frac{2d-1}{2(d-m)-1}z_0(z_{m})^{\frac{k-3}{2}-1}.
\label{final}
\end{eqnarray}
Next step goes in the same way as this step. In sum, each step changes one $z_m$ in the last factor of (\ref{final}) into $z_0$ and 
makes  the factor $\frac{1}{2(d-m)-1}$.
Therefore, we obtain the following equality.
\begin{eqnarray}
&&\frac{1}{2^{\frac{k-3}{2}}}\tilde{L}_{0}^{k,k}(e^x)\int_{-\infty}^{x}dx_1 \tilde{L}_{1}^{k,k}(e^{x_1})\int_{-\infty}^{x_1}dx_2 \tilde{L}_{2}^{k,k}(e^{x_{2}})\cdots
\int_{-\infty}^{x_{\frac{k-5}{2}}}dx_{\frac{k-3}{2}} \tilde{L}_{\frac{k-3}{2}}^{k,k}(e^{x_{\frac{k-3}{2}}})G_{o}^{k}(x_{\frac{k-3}{2}})\no\\
&=&\sum_{d=1}^{\infty}e^{(d-\frac{1}{2})x}
\sum_{m=0}^{d-1}\frac{1}{(2\pi\sqrt{-1})^{m+1}}\oint_{C_{m}}\frac{dz_{m}}{(z_{m})^k}
\oint_{C_{m-1}}\frac{dz_{m-1}}{(z_{m-1})^k}
\cdots\oint_{C_{0}}\frac{dz_{0}}{(z_{0})^k}f_{2d-2m-1}(z_0)\bigl(\prod_{i=1}^{m}e^k(z_{i-1},z_i)\bigr)\times\no\\
&&\times\frac{1}{kz_0(\frac{2}{2(d-m)-1}z_0+z_0-z_{1})}\bigl(\prod_{i=1}^{m-1}\frac{1}{kz_i(2z_i-z_{i-1}-z_{i+1})}\bigr)
\frac{2d-1}{(2(d-m)-1)^{\frac{k-1}{2}}}(z_0)^{\frac{k-3}{2}}.
\label{last}
\end{eqnarray}
We can easily see that the summands corresponding to $m>0$ vanish from power counting of the variable $z_0$, and we 
have, 
\begin{eqnarray}
&&\frac{1}{2^{\frac{k-3}{2}}}\tilde{L}_{0}^{k,k}(e^x)\int_{-\infty}^{x}dx_1 \tilde{L}_{1}^{k,k}(e^{x_1})\int_{-\infty}^{x_1}dx_2 \tilde{L}_{2}^{k,k}(e^{x_{2}})\cdots
\int_{-\infty}^{x_{\frac{k-5}{2}}}dx_{\frac{k-3}{2}} \tilde{L}_{\frac{k-3}{2}}^{k,k}(e^{x_{\frac{k-3}{2}}})G_{o}^{k}(x_{\frac{k-3}{2}})\no\\
&=&\sum_{d=1}^{\infty}e^{(d-\frac{1}{2})x}\frac{1}{(2\pi\sqrt{-1})}\oint_{C_{0}}\frac{dz_{0}}{(z_{0})^k}
f_{2d-1}(z_0)\frac{1}{(2d-1)^{\frac{k-3}{2}}}(z_0)^{\frac{k-3}{2}}=\tau^{k}(x).
\end{eqnarray} 
This equality completes the proof. $\Box$

$w_{disk}^{k,k}({\cal O}_{h^{\frac{k-3}{2}}})_{2d-1}$ is the B-model analogue of the open Gromov-Witten invariant 
$\langle{\cal O}_{h^{\frac{k-3}{2}}}\rangle_{disk,2d-1}$ and it is translated into $\langle{\cal O}_{h^{\frac{k-3}{2}}}\rangle_{disk,2d-1}$ through the mirror transformation. For general $N$ and $k$, we propose the
generalized mirror transformation for the open Gromov-Witten invariants, which is a straightforward generalization of the one 
for closed Gromov-Witten invariants \cite{gene}.
To write down the formula, we introduce the partition $\sigma_{d}$ of positive integer $d$,
\begin{eqnarray}
P_d:=\{\sigma_{d}=(d_1,d_2,\cdots,d_{l(\sigma_{d})})\;|\;1\leq d_1\leq d_2\leq\cdots\leq d_{l(\sigma_{d})},\;
d_1+ d_2+\cdots+d_{l(\sigma_{d})}=d\},
\label{partition}
\end{eqnarray}
and the symmetric factor:
\begin{eqnarray}
S(\sigma_{d}):=\prod_{i=1}^{d}\frac{1}{(\mbox{mul}(i,\sigma_d))!},
\label{symmetric}
\end{eqnarray}
where $\mbox{mul}(i,\sigma_d)$ is multiplicity of $i$ in $\sigma_{d}$.
With this set-up, the formula of the generalized mirror transformation of disk amplitude is given as follows.
\begin{conj}
{\bf (Generalized Mirror Transformation for Open Gromov-Witten invariants)}
\begin{eqnarray}
&&w_{disk}^{N,k}({\cal O}_{h^a})_{2d-1}=\no\\
&&=\langle{\cal O}_{h^{a}}\rangle_{disk,2d-1}
+\sum_{f=1}^{d-1}\sum_{\sigma_f\in P_f}S(\sigma_{f})
\langle{\cal O}_{h^{a}}\prod_{j=1}^{l(\sigma_{f})}{\cal O}_{h^{1+(k-N)f_j}}\rangle_{disk,2d-2f-1}
\prod_{j=1}^{l(\sigma_{f})}\frac{w({\cal O}_{h^{N-3-(k-N)f_j}}{\cal O}_{h^0})_{0,f_j}}{k},\no\\
\label{gmt}
\end{eqnarray}
where $w({\cal O}_{h^{N-3-(k-N)d}}{\cal O}_{h^0})_{0,d}$ is given by (\ref{wint}).
\end{conj}
In Section 5, we prove (\ref{gmt}) up to $2d-1=5$. 
This formula includes multi-point open Gromov-Witten invariant $\langle\prod_{i=1}^{n}{\cal O}_{h^{a_i}}\rangle_{disk,2d-1}$,
which is obtained by integrating out $e(F_{2d-1})\prod_{i=1}^{n}ev_{i}^{*}(h^{a_i})$ on $\widetilde{M}_{D,n}(CP^{N-1}/CP_{R}^{N-1},2d-1)$.
In this case, all the $n$ marked points are located inside the disk. For lower 
degrees, we can compute $\langle\prod_{i=1}^{n}{\cal O}_{h^{a_i}}\rangle_{disk,2d-1}$ by localization computation. In Section 5, we present 
explicit formulas to compute it up to $2d-1=5$. 

If $N=k$, we can use the K\"ahler equation for the open Gromov-Witten invariants:
\begin{equation}
\langle{\cal O}_{h^{a}}({\cal O}_{h})^{l(\sigma_{f})}\rangle_{disk,2d-2f-1}=
\langle{\cal O}_{h^{a}}\rangle_{disk,2d-2f-1}\cdot\bigl(\frac{2d-2f-1}{2}\bigr)^{l(\sigma_{f})},
\end{equation}
that can be easily verified by localization computation.
Then (\ref{gmt}) implies,
\begin{eqnarray}
\sum_{d=1}^{\infty}\langle{\cal O}_{h^{\frac{k-3}{2}}}\rangle_{disk,2d-1}e^{(d-\frac{1}{2})t(x)}=\sum_{d=1}^{\infty}w_{disk}^{k,k}({\cal O}_{h^{\frac{k-3}{2}}})_{2d-1}e^{(d-\frac{1}{2})x},
\label{CYmt}
\end{eqnarray}
where
\begin{equation}
t(x):=x+\sum_{d=1}^{\infty}\frac{w({\cal O}_{h^{k-3}}{\cal O}_{h^0})_{0,d}}{k}e^{dx}=x+
\sum_{d=1}^{\infty}\frac{\tilde{L}^{k,k,d}_{1}}{d}e^{dx}.
\end{equation}
(\ref{CYmt}) and Theorem 1 lead us to Conjecture 1, i.e., the B-model computation of the generating function of $\langle{\cal O}_{h^{\frac{k-3}{2}}}
\rangle_{disk,2d-1}$ of the Calabi-Yau hypersurface $M_{k}^{k}$.

Obviously, $ \langle{\cal O}_{h^{a}}\prod_{j=1}^{l(\sigma_{f})}{\cal O}_{h^{1+(k-N)f_i}}\rangle_{2d-2f-1}=0$ 
if $1+(k-N)f_i<0$. Hence we obtain,   
\begin{cor}
(i) If $N-k\geq 2$ and $\frac{N-3+(2d-1)(N-k)}{2}$ is integer and $k$ is odd, we have,
\begin{equation}
 \langle{\cal O}_{h^{\frac{N-3+(2d-1)(N-k)}{2}}}
\rangle_{disk,2d-1}=w_{disk}^{N,k}({\cal O}_{h^{\frac{N-3+(2d-1)(N-k)}{2}}})_{2d-1}.\\
\label{fano2}
\end{equation}
(ii) If $N-k=1$ and $\frac{N-3+(2d-1)(N-k)}{2}$ is integer and $k$ is odd, we have the following equality.
\begin{eqnarray}
w_{disk}^{N,k}({\cal O}_{h^{\frac{N-3+(2d-1)(N-k)}{2}}})_{2d-1}
&=&\sum_{j=0}^{d-1}\langle{\cal O}_{h^{\frac{N-3+(2d-1)(N-k)}{2}}}({\cal O}_{h^{0}})^{j}\rangle_{disk,2d-2j-1}
\cdot(\frac{w({\cal O}_{h^{N-2}}{\cal O}_{h^0})_{0,1}}{k})^{j}\frac{1}{j!}\no\\
&=&\sum_{j=0}^{d-1}\langle{\cal O}_{h^{\frac{N-3+(2d-1)(N-k)}{2}}}({\cal O}_{h^{0}})^{j}\rangle_{disk,2d-2j-1}
\cdot(k!)^{j}\frac{1}{j!}.
\label{fano1}
\end{eqnarray}
Here, we used the following equality:
\begin{equation}
\frac{w({\cal O}_{h^{N-2}}{\cal O}_{h^0})_{0,1}}{k}=\frac{1}{k}\frac{1}{(2\pi\sqrt{-1})^2}
\oint_{C_{(0)}}\frac{dz_{0}}{(z_{0})^{k+1}}\oint_{C_{(0)}}\frac{dz_{1}}{(z_{1})^{k+1}}e^{k}(z_{0},z_{1})\cdot(z_{0})^{k-1}=k!.
\end{equation}
\end{cor}
From (\ref{fano2}), we can see that  $\langle{\cal O}_{h^{\frac{N-3+(2d-1)(N-k)}{2}}}
\rangle_{disk,2d-1}$ is evaluated via the formula (\ref{def3}). 

\begin{Rem}
To compute  $\langle{\cal O}_{h^{\frac{N-3+(2d-1)(N-k)}{2}}}\rangle_{disk,2d-1}$ without the localization theorem in the cases of $N-k=1$ and  $N-k<0$, we 
need another equation for the open Gromov-Witten invariants to evaluate multi-point open Gromov-Witten invariants from the 1-point open Gromov-Witten invariants (it is expected to play the same role as the associativity equation in the closed string case).  
It is quite non-trivial and we leave pursuit for these topics for future works. 
\end{Rem}
\section{Numerical Data}\label{DATA}
In this section, we present the numerical data predicted by Conjectures in Section \ref{B} and Section \ref{OVSC}. 
We display the disk invariants for several hypersurfaces and complete intersections. 
The formulas for complete intersections can be obtained by straightforward extension of those of hypersurfaces. 
As is shown below, we certainly obtain integral invariants after re-summation by the (modified) multiple covering formula. 
Integral property of these invariants supports validity of our formalisms.

\subsection{Multiple Covering Formula for Calabi-Yau Hypersurfaces}
The multiple covering formula of the open Gromov-Witten invariants for Calabi-Yau $3$-folds 
was first given in the context of BPS states counting \cite{OV}. 
The formula relevant to the case of $1$-point function has the following form: 
Let $q=e^t$, where $t$ is K\"ahler moduli parameter in the A-model. 
Then, 

\begin{equation}
\sum_{d} {\langle{\cal O}_{h}\rangle}_{disk,2d-1} q^{\frac{2d-1}{2}}=\sum_{l,d}\frac{(2d-1)n_{2d-1}}{2l-1}q^{\frac{(2d-1)(2l-1)}{2}}. 
\end{equation}
Here, summation with respect to $d$ and $l$ are taken for all positive integers. 
$n_{2d-1}$ is the disk invariant of Calabi-Yau 3-fold of degree $2d-1$ and 
is conjectured to be an integer. 

Interestingly, 
we observe that the multiple covering formula is modified for the higher dimensional cases. 
For $2D+3$ ($D=1$, $2$, $3$, $...$) dimensional cases, 
we conjecture that the multiple covering formula have the following form: 
\begin{equation}
\sum_{d} {\langle{\cal O}_{h^{a}}\rangle}_{disk,2d-1} q^{\frac{2d-1}{2}}=\sum_{l,d}\frac{(-1)^{lD}n_{2d-1}}{2l-1}q^{\frac{(2d-1)(2l-1)}{2}}. 
\label{hmul}
\end{equation}
The sign $(-1)^{lD}$ in (\ref{hmul}) was found from requirement that $n_{2d-1}$ should be an integer. 
In the following, we present the integral disk invariants obtained by this formula. 
\subsection{CY $5$-fold}
\begin{center}
\begin{tabular}{|c|l|}
\hline
\multicolumn{2}{|c|}{$N=k=7$}\\ \hline
d & disk invariant \\ \hline
1 & 210 \\ [-2pt] 
3 & 20238540 \\ [-2pt]  
5 & 7164717071610 \\ [-2pt]  
7 & 3323817979294765050 \\ [-2pt] 
9 & 1753815102150400195049220 \\ [-2pt] 
11 & 997489630646125057277538604350 \\ [-2pt] 
13 & 595872023331262091783971492294372080 \\ [-2pt] 
15 & 368492591032299305435217331096538887611570 \\ [-2pt] 
17 & 233812266607099277659029702498147934247411056500 \\ [-6pt] 
$\vdots$ &  \\ \hline
\end{tabular}
\\
\vspace{2em}
\hspace{2em}
\begin{tabular}{|c|l|}
\hline
\multicolumn{2}{|c|}{$N=8$, $k_1=3$, $k_2=5$}\\ \hline
d & disk invariant \\ \hline
1 & 90 \\ [-2pt] 
3 & 819180 \\ [-2pt]  
5 & 29501336250 \\ [-2pt]  
7 & 1385913817885770 \\ [-2pt] 
9 & 74029249898896159800 \\ [-2pt] 
11 & 4262390623679148176523210 \\ [-2pt] 
13 & 257780912852264319358790011920 \\ [-2pt] 
15 & 16140190288299711750115529913146250 \\ [-2pt] 
17 & 1036944919162381766672552243000832999630 \\ [-6pt] 
$\vdots$ &  \\ \hline
\end{tabular}
\hspace{2em}
\begin{tabular}{|c|l|}
\hline
\multicolumn{2}{|c|}{$N=9$, $k_1=k_2=k_3=3$}\\ \hline
d & disk invariant \\ \hline
1 & 54 \\ [-2pt] 
3 & 106920 \\ [-2pt]  
5 & 894937302 \\ [-2pt]  
7 & 9719866853226 \\ [-2pt] 
9 & 119971275497509464 \\ [-2pt] 
11 & 1596025323666255425058 \\ [-2pt] 
13 & 22302433567318014407088792 \\ [-2pt] 
15 & 322654273568938603304257219182 \\ [-2pt] 
17 & 4789894590347957070773256970063362 \\ [-6pt] 
$\vdots$ &  \\ \hline
\end{tabular}
\end{center}
\subsection{CY $7$-fold}
\begin{center}
\begin{tabular}{|c|l|}
\hline
\multicolumn{2}{|c|}{$N=k=9$}\\ \hline
d & disk invariant \\ \hline
1 & 1890 \\ [-2pt] 
3 & 94563624960 \\ [-2pt] 
5 & 16211885196706741080 \\ [-2pt] 
7 & 3725578314401332796504317080 \\ [-2pt] 
9 & 980933196421736367520298463232432740 \\ [-2pt] 
11 & 279404400332958478341372061714624713007501620 \\ [-2pt] 
13 & 83765316559743112613635047261440115474721993367743140 \\ [-2pt] 
15 & 26032698632451236276838009397731526684304250664987130967951660 \\ [-2pt] 
17 & 8308974360738906157282140936951112739423372574571446836487783983771670 \\ [-6pt] 
$\vdots$ & \\ \hline
\end{tabular}
\\
\vspace{2em}
\begin{tabular}{|c|l|}
\hline
\multicolumn{2}{|c|}{$N=10$, $k_1=3$, $k_2=7$}\\ \hline
d & disk invariant \\ \hline
1 & 630 \\ [-2pt] 
3 & 1732513860 \\ [-2pt] 
5 & 16948700697790260 \\ [-2pt] 
7 & 222352925468971110069060 \\ [-2pt] 
9 & 3343542378535393312371806938050 \\ [-2pt] 
11 & 54403982118940628184121440657132431130 \\ [-2pt] 
13 & 931897834382491611934844844028931350370207490 \\ [-2pt] 
15 & 16549508759112266329543165383598816943011383280906290 \\ [-2pt] 
17 & 301867734062086886863292457315822598044420424688298906794790 \\ [-6pt] 
$\vdots$ &  \\ \hline
\end{tabular}
\\
\vspace{2em}
\begin{tabular}{|c|l|}
\hline
\multicolumn{2}{|c|}{$N=10$, $k_1=k_2=5$}\\ \hline
d & disk invariant \\ \hline
1 & 450 \\ [-2pt]
3 & 541083600 \\ [-2pt]
5 & 2323075852800000 \\ [-2pt]
7 & 13376764782851724374400 \\ [-2pt]
9 & 88290869224369416996530295000 \\ [-2pt]
11 & 630597905473674189253142188226611800 \\ [-2pt]
13 & 4741451594869506232624087439647902217296600 \\ [-2pt]
15 & 36961988089850513025172698893989264922228539275000 \\ [-2pt]
17 & 295950203480857846231677257507873272845493549136884959900 \\ [-6pt]
$\vdots$ &  \\ \hline
\end{tabular}
\end{center}
\subsection{CY $9$-fold}
\begin{center}
\begin{tabular}{|c|l|}
\hline
\multicolumn{2}{|c|}{$N=k=11$}\\ \hline
d & disk invariant \\ \hline
1 & 20790 \\ [-2pt] 
3 & 739689094281060 \\ [-2pt] 
5 & 92349241505201808072653400 \\ [-2pt] 
7 & 15662627763584777441409921423519150000 \\ [-2pt] 
9 & 3057781689155349311391320055260969435477344733880 \\ [-2pt] 
11 & 647228065882996100723626284929495089845800978296255510611900 \\ [-2pt] 
13 & 144373392926266851346933782476456863022618317345113666089276017766556180 \\ [-2pt] 
15 & 33410445239945586400783564864831299361460754601690634625359260018895998949763278400 \\ [-2pt] 
17 & 7944798193890486076218127410555172778932262204508240270452127481115279486189452580187157349820 \\ [-6pt] 
$\vdots$ &  \\ \hline
\end{tabular}
\\
\vspace{2em}
\begin{tabular}{|c|l|}
\hline
\multicolumn{2}{|c|}{$N=12$, $k_1=3$, $k_2=9$}\\ \hline
d & disk invariant \\ \hline
1 & 5670 \\ [-2pt] 
3 & 7190951406660 \\ [-2pt]  
5 & 32771808534815587472760 \\ [-2pt]  
7 & 203196564019239912556407749039760 \\ [-2pt] 
9 & 1451156382550401076913171359345180261442280 \\ [-2pt] 
11 & 11239991067300623274937381702339644645299208276030540 \\ [-2pt] 
13 & 91766796783481120695413660015328806345445819073092114838184100 \\ [-2pt] 
15 & 777374205495124523505951877224678810574278189880401144822687459772831920 \\ [-2pt] 
17 & 6767415192834983648586789918872473794816814205793015494848356268766369435917680660 \\ [-6pt] 
$\vdots$ &  \\ \hline
\end{tabular}
\\
\vspace{2em}
\begin{tabular}{|c|l|}
\hline
\multicolumn{2}{|c|}{$N=12$, $k_1=5$, $k_2=7$}\\ \hline
d & disk invariant \\ \hline
1 & 3150 \\ [-2pt] 
3 & 979005705300 \\ [-2pt]  
5 & 1096990088649650685000 \\ [-2pt]  
7 & 1672634922992392369191998947200 \\ [-2pt] 
9 & 2937773983965272566082324281789596936000 \\ [-2pt] 
11 & 5596418777758382921429206857192171069228643754100 \\ [-2pt] 
13 & 11237841476030778718599291797560457675263594779179678570700 \\ [-2pt] 
15 & 23414742446591074932980593239147196095984682847992899044385164745000 \\ [-2pt] 
17 & 50136079032189916998378367844077844900721498505814918103036703603401027454800 \\ [-6pt] 
$\vdots$ &  \\ \hline
\end{tabular}
\end{center}
\subsection{CY $11$-fold}
\begin{center}
\begin{tabular}{|c|l|}
\hline
\multicolumn{2}{|c|}{$N=k=13$}\\ \hline
d & disk invariant \\ \hline
1 & 270270 \\ [-2pt] 
3 & 9630787776863673420 \\ [-2pt]  
5 & 1259056659533544456991412149863720 \\ [-2pt]  
7 & 225580159504495415590622751696044661618421816800 \\ [-2pt] 
9 & 46657333666806812735845639185680874962204077273118803505703120 \\ [-2pt] 
11 & 10476928443077114442692287150154959225437315286976833827059990771981906214300 \\ [-2pt] 
13 &248118594108238306151718332183554970781499917826353464617784863271086518825323-\\[-2pt] 
& -3441903339580\\ [-2pt] 
15 &609904678687013223891711069528390016779635460997454892925440741188365241178962-\\[-2pt] 
& -769648409087639739094525960\\ [-2pt] 
17 &154104919821094255782922784038386351047258709232740417178312240458745268535691-\\[-2pt] 
&-125346421511857117632677043395862969090200 \\ [-6pt] 
$\vdots$ &  \\ \hline
\end{tabular}
\end{center}
\subsection{Fano Hypersurfaces }
\begin{center}
\begin{tabular}{|c|l|l|}
\hline
\multicolumn{3}{|c|}{$N=5$, $k=3$}\\ \hline
d & $\langle{\cal O}_{h^{\frac{N-3+(N-k)d}{2}}}\rangle_{disk,d}$ & 
$w_{disk}^{5,3}({\cal O}_{h^{\frac{N-3+(N-k)d}{2}}})_{d}$\\ \hline
1 & 6 &6\\ \hline
\end{tabular}
\\
\vspace{2em}
\begin{tabular}{|c|l|l|}
\hline
\multicolumn{3}{|c|}{$N=6$, $k=5$}\\ \hline
d & $\langle{\cal O}_{h^{\frac{N-3+(N-k)d}{2}}}\rangle_{disk,d}$  & $w_{disk}^{6,5}({\cal O}_{h^{\frac{N-3+(N-k)d}{2}}})_{d}$\\ \hline
1 & 30 & 30　\\ [-2pt] 
3 & 6300 & 8100　　\\ [-2pt]  
5 & 1198800 & 1216800　　\\ 
\hline
\end{tabular}
\hspace{2em}
\begin{tabular}{|c|l|l|}
\hline
\multicolumn{3}{|c|}{$N=7$, $k=5$}\\ \hline
d & $\langle{\cal O}_{h^{\frac{N-3+(N-k)d}{2}}}\rangle_{disk,d}$ & $w_{disk}^{7,5}({\cal O}_{h^{\frac{N-3+(N-k)d}{2}}})_{d}$\\ \hline
1 & 30 & 30\\ [-2pt] 
3 & 1200 & 1200  \\  
\hline
\end{tabular}
\\
\vspace{2em}
\begin{tabular}{|c|l|l|}
\hline
\multicolumn{3}{|c|}{$N=8$, $k=7$}\\ \hline
d & $\langle{\cal O}_{h^{\frac{N-3+(N-k)d}{2}}}\rangle_{disk,d}$ & $w_{disk}^{8,7}({\cal O}_{h^{\frac{N-3+(N-k)d}{2}}})_{d}$ \\ \hline
1 & 210 &　210 \\ [-2pt] 
3 & 14852880 &　15382080 \\ [-2pt]  
5 & 1082061907920 &　1091604442320 \\ [-2pt]  
7 & 20924080987824000 &　21063211139376000  \\ 
\hline
\end{tabular}
\hspace{2em}
\begin{tabular}{|c|l|l|}
\hline
\multicolumn{3}{|c|}{$N=9$, $k=7$}\\ \hline
d & $\langle{\cal O}_{h^{\frac{N-3+(N-k)d}{2}}}\rangle_{disk,d}$ & $w_{disk}^{9,7}({\cal O}_{h^{\frac{N-3+(N-k)d}{2}}})_{d}$\\ \hline
1 & 210 & 210\\ [-2pt] 
3 & 4051320 & 4051320\\ 
 \hline
\end{tabular}
\vspace{1em}
\end{center}
In the case of Fano hypersurfaces, $\langle{\cal O}_{h^{\frac{N-3+(N-k)d}{2}}}\rangle_{disk,d}$ is evaluated by localization computation
(see the formulas (\ref{ogres}) in Section 5). If $N-k=1$, $\langle{\cal O}_{h^{\frac{N-3+(N-k)d}{2}}}\rangle_{disk,d}$ is different from $w_{disk}({\cal O}_{h^{\frac{N-3+(N-k)d}{2}}})_{d}$ because of the generalized mirror transformation in (\ref{fano1}). Let us see the case of $N=8,k=7,d=7
$ in detail. In this case, $\langle{\cal O}_{h^{6}}{\cal O}_{h^{0}}\rangle_{disk,5}$, $\langle{\cal O}_{h^{6}}({\cal O}_{h^{0}})^2\rangle_{disk,3}$ and $\langle{\cal O}_{h^{6}}({\cal O}_{h^{0}})^3\rangle_{disk,1}$ are evaluated from the formulas (\ref{ogres})
in Section 5 as follows. 
\begin{eqnarray}
\langle{\cal O}_{h^{6}}{\cal O}_{h^{0}}\rangle_{disk,5}=27605188800,\;
\langle{\cal O}_{h^{6}}({\cal O}_{h^{0}})^2\rangle_{disk,3}=-44100, \langle{\cal O}_{h^{6}}({\cal O}_{h^{0}})^3\rangle_{disk,1}=105/4.
\label{example87}
\end{eqnarray}
With these data, we can confirm (\ref{fano1}) in this case explicitly.
\begin{equation}
w_{disk}^{8,7}({\cal O}_{h^{6}})_{d}=21063211139376000=20924080987824000+27605188800\cdot(7!)+(-44100)\cdot(7!)^2/2+(105/4)\cdot(7!)^3/6.
\end{equation}
Lastly, we present the numerical data for a general type hypersurface with $N=8,k=9$. In this case, we have three 1 point disk Gromov-Witten 
invariants.
\begin{eqnarray}
\langle{\cal O}_{h^{2}}\rangle_{disk,1}=1890,\;\langle{\cal O}_{h^{1}}\rangle_{disk,3}=58381461390,\;
\langle{\cal O}_{h^{0}}\rangle_{disk,5}=41731576876146796884/25.
\end{eqnarray}
The corresponding open virtual structure constants are given as follows.
\begin{eqnarray}
w_{disk}^{8,9}({\cal O}_{h^{2}})_{1}=1890,\;w_{disk}^{8,9}({\cal O}_{h^{1}})_{3}=90642729450,\;w_{disk}^{8,9}({\cal O}_{h^{0}})_{5}=276177175032776063634/25. 
\end{eqnarray}
To confirm (\ref{gmt}) numerically, we prepare the following data evaluated from the formulas (\ref{ogres}) in Section 5 and (\ref{wint}).
\begin{eqnarray}
&&\langle {\cal O}_{h^{1}}{\cal O}_{h^{2}} \rangle_{disk,1}=945,\;
\langle{\cal O}_{h^{0}}{\cal O}_{h^{3}}\rangle_{disk,1}=945,\;
\langle{\cal O}_{h^{0}}({\cal O}_{h^{2}})^2\rangle_{disk,1}=945/2,\;
\langle{\cal O}_{h^{0}}{\cal O}_{h^{2}}\rangle_{disk,3}=33973546005,\no\\
&&\frac{w({\cal O}_{h^{0}}{\cal O}_{h^{4}})_{0,1}}{9}=34138908,\;
\frac{w({\cal O}_{h^{0}}{\cal O}_{h^{3}})_{0,2}}{9}=8404934443598718,\;
\end{eqnarray}
Then (\ref{gmt}) in this case reduces to the following equalities:
\begin{eqnarray}
&&w_{disk}^{8,9}({\cal O}_{h^{2}})_{1}=1890=1890,\no\\
&&w_{disk}^{8,9}({\cal O}_{h^{1}})_{3}=90642729450=58381461390+945\cdot34138908,\no\\ 
&&w_{disk}^{8,9}({\cal O}_{h^{0}})_{5}=276177175032776063634/25=41731576876146796884/25+33973546005\cdot 34138908+\no\\
&&\hspace{7em}945\cdot 8404934443598718+(945/2)\cdot(34138908)^2/2.
\end{eqnarray} 
\section{Proof of the Generalized Mirror Transformation up to the $d=5$ Case}\label{PROOF}
In this section, we first present the formulas that compute the open Gromov-Witten invariants of $M_{N}^{k}$ up to 
$d=5$. Next, we prove the conjectures in Section \ref{OVSC} up to $d=5$. For this purpose, 
we introduce the following polynomial in $x$ and $y$:
\begin{eqnarray}
h_{a}(x,y)&:=&\frac{x^a-y^a}{x-y},
\end{eqnarray}
where $a$ is a non-negative integer.
Our formulas to compute the open Gromov-Witten invariants up to $\displaystyle{d=5}$ are given as follows.
\begin{prop}
The A-model amplitude $\langle\prod_{i=1}^{n}{\cal O}_{h^{a_i}}\rangle_{disk,d}$ up to the $d=5$ case 
is given by sum of the following residue integrals 
associated with graphs in Figure 3.   
\begin{eqnarray}
\langle\prod_{i=1}^{n}{\cal O}_{h^{a_i}}\rangle_{disk,1}&=&\frac{1}{(2\pi\sqrt{-1})}\oint_{C_{0}}\frac{dz_0}{(z_{0})^{N}}
f_{1}(z_{0})\cdot2z_0\cdot\prod_{i=1}^{n}(\frac{(z_0)^{a_i-1}}{2}),\no\\
\langle\prod_{i=1}^{n}{\cal O}_{h^{a_i}}\rangle_{disk,3}&=&\frac{1}{(2\pi\sqrt{-1})}\oint_{C_{0}}\frac{dz_0}{(z_{0})^{N}}
f_{3}(z_{0})\cdot\frac{2z_0}{3}\cdot\prod_{i=1}^{n}(\frac{3(z_0)^{a_i-1}}{2})\no\\
&+&\frac{1}{(2\pi\sqrt{-1})^2}\oint_{C_{0}}\frac{dz_0}{(z_{0})^{N}}\oint_{C_{1}}\frac{dz_1}{(z_{1})^{N}}
f_{1}(z_{0})e^{k}(z_0,z_1)\frac{z_1-z_0}{kz_0(3z_0-z_1)}
\cdot\prod_{i=1}^{n}(\frac{(z_0)^{a_i-1}}{2}+h_{a_i}(z_0,z_1)), \no\\
\langle\prod_{i=1}^{n}{\cal O}_{h^{a_i}}\rangle_{disk,5}&=&\frac{1}{(2\pi\sqrt{-1})}\oint_{C_{0}}\frac{dz_0}{(z_{0})^{N}}
f_{5}(z_{0})\cdot\frac{2z_0}{5}\cdot\prod_{i=1}^{n}(\frac{5(z_0)^{a_i-1}}{2})\no\\
&+&\frac{1}{(2\pi\sqrt{-1})^2}\oint_{C_{0}}\frac{dz_0}{(z_{0})^{N}}\oint_{C_{1}}\frac{dz_1}{(z_{1})^{N}}
f_{3}(z_{0})e^{k}(z_0,z_1)\frac{z_1-z_0}{kz_0(\frac{5}{3}z_0-z_1)}
\cdot \prod_{i=1}^{n}(\frac{3(z_0)^{a_i-1}}{2}+h_{a_i}(z_0,z_1))\no\\
&+&\frac{1}{(2\pi\sqrt{-1})^3}\oint_{C_{0}}\frac{dz_0}{(z_{0})^{N}}\oint_{C_{1}}\frac{dz_1}{(z_{1})^{N}}
\oint_{C_{2}}\frac{dz_2}{(z_{2})^{N}}
f_{1}(z_{0})e^{k}(z_0,z_1)e^{k}(z_1,z_2)\times\no\\
&&\times\frac{z_2-z_1}{kz_0(3z_0-z_1)kz_1(2z_1-z_0-z_2)}
\prod_{i=1}^{n}(\frac{(z_0)^{a_i-1}}{2}+h_{a_i}(z_0,z_1)+h_{a_i}(z_1,z_2))\no\\
&+&\frac{1}{2}\frac{1}{(2\pi\sqrt{-1})^3}\oint_{C_{0}}\frac{dz_0}{(z_{0})^{N}}\oint_{C_{1}}\frac{dz_1}{(z_{1})^{N}}
\oint_{C_{2}}\frac{dz_2}{(z_{2})^{N}}
f_{1}(z_{0})e^{k}(z_0,z_1)e^{k}(z_0,z_2)\frac{1}{(kz_0)^2(2z_0)}\times\no\\
&&\times\prod_{i=1}^{n}(\frac{(z_0)^{a_i-1}}{2}+h_{a_i}(z_0,z_1)+h_{a_i}(z_0,z_2)).
\label{ogres}
\end{eqnarray}
\end{prop}
\begin{figure}[h]
      \epsfxsize=12cm
     \centerline{\epsfbox{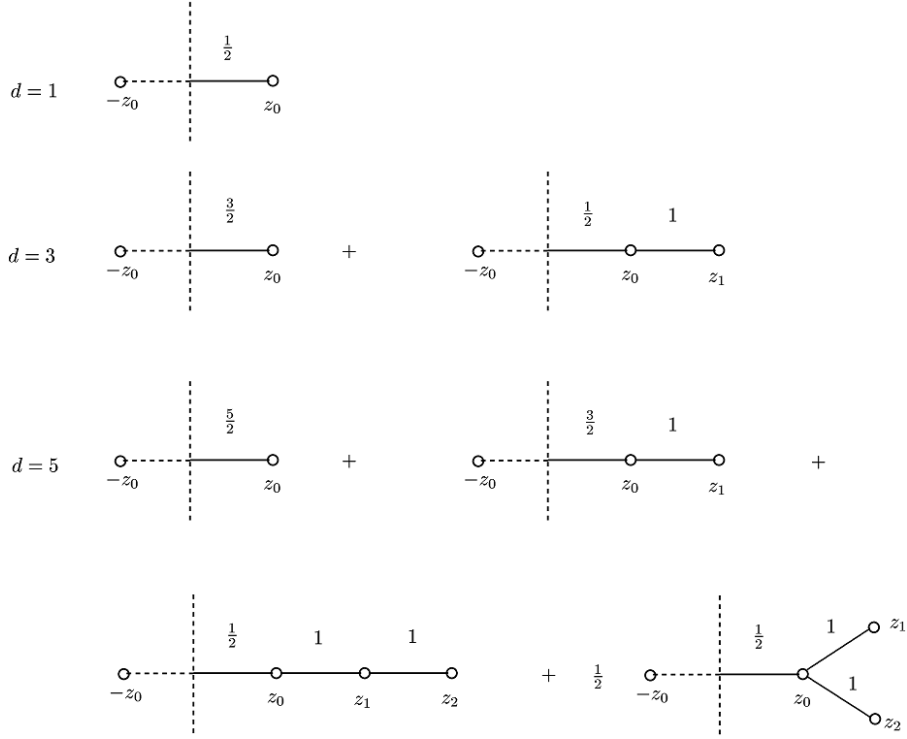}}
    \caption{\bf The graphs that contribute to A-model amplitudes }
\label{gA}
\end{figure}
In the above formulas, we take the residue integrals in ascending order of the subscript $i$ of $z_i$. $\frac{1}{2\sqrt{-1}}\oint_{C_i}dz_i$
means that we take the residues at $z_{i}=0,\frac{z_{i-1}+z_{i+1}}{2}$ (resp. $z_i=0$) if the integrand contains the factor
$\frac{1}{2z_i-z_{i-1}-z_{i+1}}$ (resp. otherwise).
This proposition follows from the localization computation applied to the open Gromov-Witten invariants \cite{kont} \cite{walcher}, the non-equivariant limit
used in \cite{pmth} and the same trick as the one that reduces (\ref{odef3}) to (\ref{def3}).
\begin{theorem}
If $d\leq 5$, (\ref{gmt}) holds true. Explicitly, we have, 
\begin{eqnarray}
w_{disk}^{N,k}({\cal O}_{h^a})_{1}&=&\langle{\cal O}_{h^{a}}\rangle_{disk,1},\\
w_{disk}^{N,k}({\cal O}_{h^a})_{3}&=&\langle{\cal O}_{h^{a}}\rangle_{disk,3}+
\langle{\cal O}_{h^{a}}{\cal O}_{h^{1+(k-N)}}\rangle_{disk,1}\frac{w({\cal O}_{h^{N-3-(k-N)}}{\cal O}_{h^0})_{0,1}}{k},\\
w_{disk}^{N,k}({\cal O}_{h^a})_{5}&=&\langle{\cal O}_{h^{a}}\rangle_{disk,5}+
\langle{\cal O}_{h^{a}}{\cal O}_{h^{1+(k-N)}}\rangle_{disk,3}\frac{w({\cal O}_{h^{N-3-(k-N)}}{\cal O}_{h^0})_{0,1}}{k}\no\\
&&+\langle{\cal O}_{h^{a}}{\cal O}_{h^{1+2(k-N)}}\rangle_{disk,1}\frac{w({\cal O}_{h^{N-3-2(k-N)}}{\cal O}_{h^0})_{0,2}}{k}\no\\ 
&&+\frac{1}{2}\langle{\cal O}_{h^{a}}({\cal O}_{h^{1+(k-N)}})^2\rangle_{disk,1}(\frac{w({\cal O}_{h^{N-3-(k-N)}}{\cal O}_{h^0})_{0,1}}{k})^2. \label{5/2}
\end{eqnarray}
\end{theorem}
{\it proof)} The case of $d=1$ is trivial. We only prove the case of $2d-1=5$. The proof for the case of 
$2d-1=3$ goes in 
the same way. First, we rewrite the second term of the R.H.S. in (\ref{5/2}), $\displaystyle{\langle{\cal O}_{h^{a}}{\cal O}_{h^{1+(k-N)}}\rangle_{disk,3}\frac{w({\cal O}_{h^{N-3-(k-N)}}{\cal O}_{h^0})_{0,1}}{k}}$, by using Proposition 1. 
\begin{eqnarray}
&&\langle{\cal O}_{h^{a}}{\cal O}_{h^{1+(k-N)}}\rangle_{disk,3}\frac{w({\cal O}_{h^{N-3-(k-N)}}{\cal O}_{h^0})_{0,1}}{k}\no\\
&=&\bigl(\frac{1}{(2\pi\sqrt{-1})}\oint_{C_{0}}\frac{dz_0}{(z_{0})^{N}}
f_{3}(z_{0})\cdot\frac{2z_0}{3}\cdot(\frac{3(z_0)^{a-1}}{2})(\frac{3(z_0)^{1+(k-N)}}{2z_0})\no\\
&&+\frac{1}{(2\pi\sqrt{-1})^2}\oint_{C_{0}}\frac{dz_0}{(z_{0})^{N}}\oint_{C_{1}}\frac{dz_1}{(z_{1})^{N}}
f_{1}(z_{0})e^{k}(z_0,z_1)\frac{z_1-z_0}{kz_0(3z_0-z_1)}(\frac{(z_0)^{a-1}}{2}+h_{a}(z_0,z_1))\times\no\\
&&\times(\frac{(z_{0})^{1+k-N}}{2z_0}+\frac{(z_0)^{1+k-N}}{z_0-z_1}+\frac{(z_1)^{1+k-N}}{z_1-z_0})\bigr)
\frac{1}{k}\bigl(\frac{1}{(2\pi\sqrt{-1})^2}\oint_{C_{0}}\frac{dz_0}{(z_{0})^{N}}\oint_{C_{1}}\frac{dz_1}{(z_{1})^{N}}
(z_0)^{N-3-(k-N)}e^{k}(z_0,z_1)\bigr)\no\\
&=&\bigl(\frac{1}{(2\pi\sqrt{-1})}\oint_{C_{0}}\frac{dz_0}{(z_{0})^{N}}
f_{3}(z_{0})\cdot\frac{2z_0}{3}\cdot(\frac{3(z_0)^{a-1}}{2})\frac{3(z_0)^{1+(k-N)}}{2z_0}\no\\
&&+\frac{1}{(2\pi\sqrt{-1})^2}\oint_{C_{0}}\frac{dz_0}{(z_{0})^{N}}\oint_{C_{1}}\frac{dz_1}{(z_{1})^{N}}
f_{1}(z_{0})e^{k}(z_0,z_1)\frac{z_1-z_0}{kz_0(3z_0-z_1)}(\frac{(z_0)^{a-1}}{2}+h_{a}(z_0,z_1))\times\no\\
&&\times(\frac{(z_{0})^{1+k-N}(3z_0-z_1)}{2z_0(z_0-z_1)}+\frac{(z_1)^{1+k-N}}{z_1-z_0})\bigr)
\frac{1}{k}\bigl(\frac{1}{(2\pi\sqrt{-1})^2}\oint_{C_{0}}\frac{dz_0}{(z_{0})^{N}}\oint_{C_{1}}\frac{dz_1}{(z_{1})^{N}}
(z_0)^{N-3-(k-N)}e^{k}(z_0,z_1)\bigr)
\end{eqnarray}
Since $(1+k-N)+(N-3-(k-N))=N-2$, we can rewrite product of two residue integrals into one residue integral by using the same 
trick as was used in the proof of Theorem 1. Hence we obtain, 
\begin{eqnarray}
&&\langle{\cal O}_{h^{a}}{\cal O}_{h^{1+(k-N)}}\rangle_{disk,3}\frac{w({\cal O}_{h^{N-3-(k-N)}}{\cal O}_{h^0})_{0,1}}{k}\no\\
&=&\frac{3}{2}\frac{1}{(2\pi\sqrt{-1})^2}\oint_{C_{0}}\frac{dz_0}{(z_{0})^{N}}\oint_{C_{1}}\frac{dz_1}{(z_{1})^{N}}
f_{3}(z_{0})\cdot\frac{2z_0}{3}\cdot(\frac{3(z_0)^{a-2}}{2})\frac{1}{kz_0}e^{k}(z_0,z_1)\no\\
&&-\frac{1}{(2\pi\sqrt{-1})^3}\oint_{C_{0}}\frac{dz_0}{(z_{0})^{N}}\oint_{C_{1}}\frac{dz_1}{(z_{1})^{N}}
\oint_{C_{2}}\frac{dz_2}{(z_{2})^{N}}
f_{1}(z_{0})\frac{1}{2z_0}\frac{1}{(kz_0)^2}e^{k}(z_0,z_1)e^{k}(z_0,z_2)(\frac{(z_0)^{a-1}}{2}+h_{a}(z_0,z_1))\no\\
&&+\frac{1}{(2\pi\sqrt{-1})^3}\oint_{C_{0}}\frac{dz_0}{(z_{0})^{N}}\oint_{C_{1}}\frac{dz_1}{(z_{1})^{N}}
\oint_{C_{2}}\frac{dz_2}{(z_{2})^{N}}
f_{1}(z_{0})\frac{1}{kz_0(3z_{0}-z_{1})}\frac{1}{kz_{1}}e^{k}(z_0,z_1)e^{k}(z_1,z_2)
(\frac{(z_0)^{a-1}}{2}+h_{a}(z_0,z_1))\no\\
&=&\frac{1}{(2\pi\sqrt{-1})^2}\oint_{C_{0}}\frac{dz_0}{(z_{0})^{N}}\oint_{C_{1}}\frac{dz_1}{(z_{1})^{N}}
f_{3}(z_{0})\cdot(\frac{3(z_0)^{a-1}}{2})\frac{1}{kz_0}e^{k}(z_0,z_1)\no\\
&&-\frac{1}{2}\frac{1}{(2\pi\sqrt{-1})^3}\oint_{C_{0}}\frac{dz_0}{(z_{0})^{N}}\oint_{C_{1}}\frac{dz_1}{(z_{1})^{N}}
\oint_{C_{2}}\frac{dz_2}{(z_{2})^{N}}
f_{1}(z_{0})\frac{1}{2z_0}\frac{1}{(kz_0)^2}e^{k}(z_0,z_1)e^{k}(z_0,z_2)\times\no\\
&&\times(2\cdot\frac{(z_0)^{a-1}}{2}+h_{a}(z_0,z_1)
+h_{a}(z_0,z_2))\no\\
&&+\frac{1}{(2\pi\sqrt{-1})^3}\oint_{C_{0}}\frac{dz_0}{(z_{0})^{N}}\oint_{C_{1}}\frac{dz_1}{(z_{1})^{N}}
\oint_{C_{2}}\frac{dz_2}{(z_{2})^{N}}
f_{1}(z_{0})\frac{1}{kz_0(3z_{0}-z_{1})}\frac{1}{kz_{1}}e^{k}(z_0,z_1)e^{k}(z_1,z_2)\times\no\\
&&\times(\frac{(z_0)^{a-1}}{2}+h_{a}(z_0,z_1)).
\end{eqnarray}
Next, we rewrite the third term, $\displaystyle{\langle{\cal O}_{h^{a}}{\cal O}_{h^{1+2(k-N)}}\rangle_{disk,1}\frac{w({\cal O}_{h^{N-3-2(k-N)}}{\cal O}_{h^0})_{0,2}}{k}}$, 
\begin{eqnarray}
&&\langle{\cal O}_{h^{a}}{\cal O}_{h^{1+2(k-N)}}\rangle_{disk,1}\frac{w({\cal O}_{h^{N-3-2(k-N)}}{\cal O}_{h^0})_{0,2}}{k}\no\\
&=&\bigl(\frac{1}{(2\pi\sqrt{-1})}\oint_{C_{0}}\frac{dz_0}{(z_{0})^{N}}
f_{1}(z_{0})\cdot2z_0\cdot(\frac{(z_0)^{a-1}}{2})\frac{(z_0)^{1+2(k-N)}}{2z_0}\bigr)\times\no\\
&&\times\frac{1}{k}\bigl(\frac{1}{(2\pi\sqrt{-1})^3}\oint_{C_{0}}\frac{dz_0}{(z_{0})^{N}}\oint_{C_{1}}\frac{dz_1}{(z_{1})^{N}}
\oint_{C_{2}}\frac{dz_2}{(z_{2})^{N}}(z_0)^{N-3-2(k-N)}
e^{k}(z_0,z_1)\frac{1}{kz_1(2z_{1}-z_{0}-z_{2})}e^{k}(z_1,z_2)\bigr)\no\\
&=&\frac{1}{(2\pi\sqrt{-1})^3}\oint_{C_{0}}\frac{dz_0}{(z_{0})^{N}}\oint_{C_{1}}\frac{dz_1}{(z_{1})^{N}}
\oint_{C_{2}}\frac{dz_2}{(z_{2})^{N}}f_{1}(z_{0})\frac{1}{kz_0}(\frac{(z_0)^{a-1}}{2})e^{k}(z_0,z_1)\frac{1}{kz_1(2z_{1}-z_{0}-z_{2})}e^{k}(z_1,z_2).\no\\
\end{eqnarray}
Finally, we rewrite the last term, $\displaystyle{\frac{1}{2}\langle{\cal O}_{h^{a}}({\cal O}_{h^{1+(k-N)}})^2\rangle_{disk,1}(\frac{w({\cal O}_{h^{N-3-(k-N)}}{\cal O}_{h^0})_{0,1}}{k})^2}$, 
\begin{eqnarray}
&&\frac{1}{2}\langle{\cal O}_{h^{a}}({\cal O}_{h^{1+(k-N)}})^2\rangle_{disk,1}(\frac{w({\cal O}_{h^{N-3-(k-N)}}{\cal O}_{h^0})_{0,1}}{k})^2
\no\\
&=&\frac{1}{2}\bigl(\frac{1}{(2\pi\sqrt{-1})}\oint_{C_{0}}\frac{dz_0}{(z_{0})^{N}}
f_{1}(z_{0})\cdot2z_0\cdot(\frac{(z_0)^{a-1}}{2})(\frac{(z_0)^{1+(k-N)}}{2z_0})^2\bigr)\times\no\\
&&\times\bigl(\frac{1}{k}\frac{1}{(2\pi\sqrt{-1})^2}\oint_{C_{0}}\frac{dz_0}{(z_{0})^{N}}\oint_{C_{1}}\frac{dz_1}{(z_{1})^{N}}
(z_0)^{N-3-(k-N)}e^{k}(z_0,z_1)\bigr)^2\no\\
&=&\frac{1}{2}\frac{1}{(2\pi\sqrt{-1})^3}\oint_{C_{0}}\frac{dz_0}{(z_{0})^{N}}\oint_{C_{1}}\frac{dz_1}{(z_{1})^{N}}
\oint_{C_{2}}\frac{dz_2}{(z_{2})^{N}}f_{1}(z_{0})\cdot(\frac{(z_0)^{a-1}}{2})\frac{1}{2z_0}\frac{1}{(kz_0)^2}
e^{k}(z_0,z_1)e^{k}(z_0,z_2).
\end{eqnarray}
\begin{figure}[h]
      \epsfxsize=12cm
     \centerline{\epsfbox{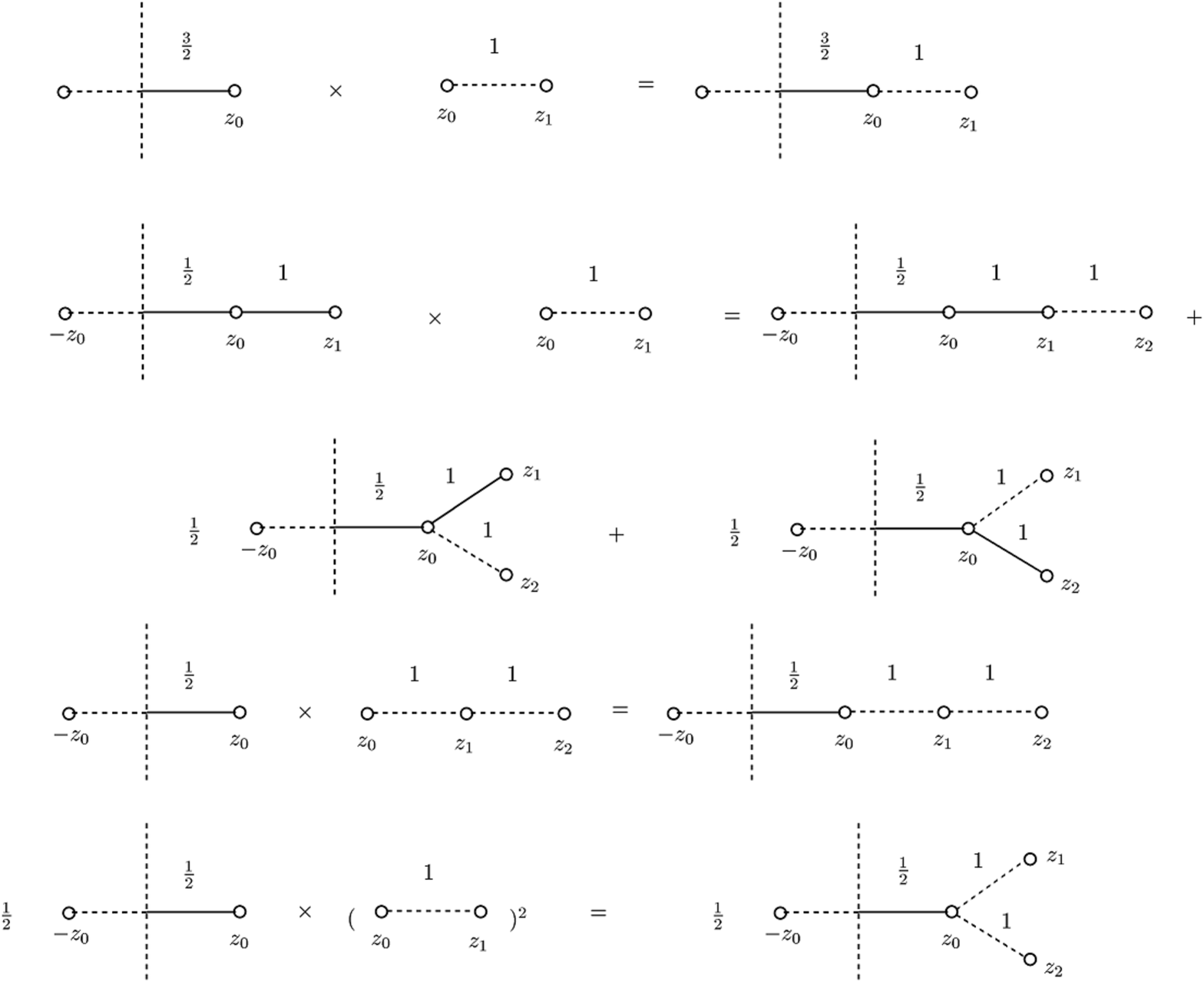}}
    \caption{\bf Graphical Representation of the Proof of Theorem 2 (1) }
\label{open-prod}
\end{figure}
\begin{figure}[h]
      \epsfxsize=12cm
     \centerline{\epsfbox{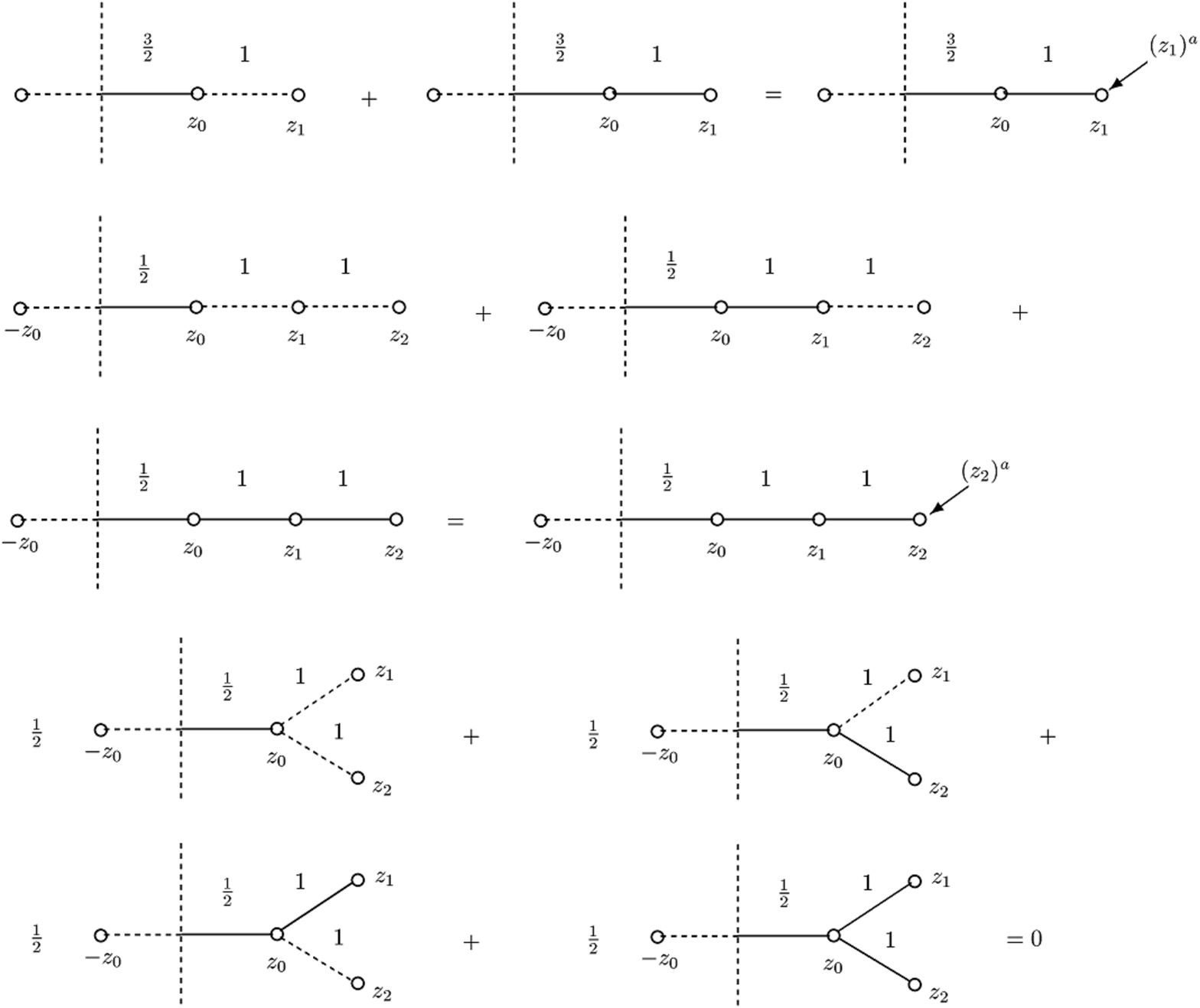}}
    \caption{\bf Graphical Representation of the Proof of Theorem 2 (2) }
\label{open-sum}
\end{figure}

On the other hand, $\langle{\cal O}_{h^{a}}\rangle_{disk,5}$ is given by,
\begin{eqnarray}
\langle{\cal O}_{h^{a}}\rangle_{disk,5}
&=&f_{5}(z_{0})\cdot\frac{2z_0}{5}\cdot(\frac{5(z_0)^{a-1}}{2})\no\\
&+&\frac{1}{(2\pi\sqrt{-1})^2}\oint_{C_{0}}\frac{dz_0}{(z_{0})^{N}}\oint_{C_{1}}\frac{dz_1}{(z_{1})^{N}}
f_{3}(z_{0})e^{k}(z_0,z_1)\frac{z_1-z_0}{kz_0(\frac{5}{3}z_0-z_1)}
(\frac{3(z_0)^{a-1}}{2}+h_{a}(z_0,z_1))\no\\
&+&\frac{1}{(2\pi\sqrt{-1})^3}\oint_{C_{0}}\frac{dz_0}{(z_{0})^{N}}\oint_{C_{1}}\frac{dz_1}{(z_{1})^{N}}
\oint_{C_{2}}\frac{dz_2}{(z_{2})^{N}}
f_{1}(z_{0})e^{k}(z_0,z_1)e^{k}(z_1,z_2)\times\no\\
&&\times\frac{z_2-z_1}{kz_0(3z_0-z_1)kz_1(2z_1-z_0-z_2)}
(\frac{(z_0)^{a-1}}{2}+h_{a}(z_0,z_1)+h_{a_i}(z_1,z_2))\no\\
&+&\frac{1}{2}\frac{1}{(2\pi\sqrt{-1})^3}\oint_{C_{0}}\frac{dz_0}{(z_{0})^{N}}\oint_{C_{1}}\frac{dz_1}{(z_{1})^{N}}
\oint_{C_{2}}\frac{dz_2}{(z_{2})^{N}}
f_{1}(z_{0})e^{k}(z_0,z_1)e^{k}(z_0,z_2)\frac{1}{(kz_0)^2(2z_0)}\times\no\\
&&\times(\frac{(z_0)^{a-1}}{2}+h_{a}(z_0,z_1)+h_{a_i}(z_0,z_2)).
\end{eqnarray}
Therefore, the theorem follows from the following elementary identities.
\begin{eqnarray}
&&(\frac{5}{3}z_0-z_1)\frac{3}{2}(z_0)^{a-1}+(z_1-z_0)\frac{3}{2}(z_0)^{a-1}+(z_1)^a-(z_0)^a
=(z_1)^a,\no\\
&&(2z_1-z_0-z_2)(\frac{1}{2}(z_0)^{a-1}+h_a(z_0,z_1))+(3z_0-z_1)\frac{1}{2}(z_0)^{a-1}+
(z_2-z_1)(\frac{1}{2}(z_0)^{a-1}+h_a(z_0,z_1)+h_a(z_1,z_2))\no\\
&&=2z_0\cdot\frac{1}{2}(z_0)^{a-1}+(z_1-z_0)h_a(z_0,z_1)+(z_2-z_1)h_a(z_1,z_2)=(z_2)^a,\no\\
&&\frac{1}{2}(z_0)^{a-1}+h_a(z_0,z_1)+h_a(z_0,z_2)-(2\cdot\frac{1}{2}(z_0)^{a-1}+h_a(z_0,z_1)+h_a(z_0,z_2))
+\frac{1}{2}(z_0)^{a-1}=0.
\end{eqnarray}

$\Box$
\begin{Rem}
The process of computation in the proof of Theorem 2 can be represented by using graphs as is shown in Figure 4 and Figure 5. 
In these figures, the parts corresponding to $\langle{\cal O}_{h^a}\prod_{i=1}^{l(\sigma_{f})}{\cal O}_{1+(k-N)f_j}\rangle_{disk,2d-2f-1}$ 
are represented by graphs with thick straight lines and the parts corresponding to $w({\cal O}_{N-3-(k-N)f_j}{\cal O}_{h^0})_{0,f_j}$ 
are represented by graphs with dashed lines. 
\end{Rem}


\newpage

\appendix
\section{Direct Integration of the Period Integral}\label{DIRECT}
In this appendix, 
we compute the period integrals of Calabi-Yau projective hypersurfaces by using (the generalized version of) the method developed in \cite{FNSS,SS}. 

Let $P$ be the defining equation of $W_k^k$ (the mirror partner of $M_k^k$): 
\begin{equation}
\{P=X_1^k+X_2^k+...+X_{k-1}^k+X_k^k-(k\psi ) X_1X_2... X_k=0\}\subset CP^{k-1}/(Z_k)^{k-2}. \label{defeqB}
\end{equation}
Now we restrict our attention to the odd degree cases, $k=2m+1$ ($m\in Z_{\geq 2}$), just as assumed throughout the paper. 
We refer \cite{GMP} for the (closed) mirror symmetry of higher dimensional Calabi-Yau manifolds. 

First, let us remember that the number of the complex structure moduli of $W_k^k$ is encoded in $h_{k-3,1}:= \mbox{rank}(H^{k-3,1}(W_k^k,Z))$. 
In our case, $h_{k-3,1}=1$ and $\psi$ in (\ref{defeqB}) represents a parameter of the complex structure. 
The relation to the parameter $x$ used in the text is $\frac{1}{(k\psi)^{k}}=e^{x}$ and we use $z=\frac{1}{(k\psi)^{k}}$ in the following. 
The large complex structure point (the mirror point of the large volume limit) is $z=0$. 
The period integrals, which play crucial roles in the B-model side of mirror symmetry, 
are defined by the integral of the holomorphic $(k-2)$-form $\Omega^{k-2,0}$ over $(k-2)$-cycles $C^{k-2}_j\in H_{k-2}(W_k^k,Z)$ ($j=0$, ..., $k-2$): 
\begin{equation}
w_j(z)=\int_{C^{k-2}_j}\Omega^{k-2,0}. \label{periodintegral}
\end{equation}
In the projective hypersurface cases, 
the holomorphic $(k-2)$-form $\Omega^{k-2,0}$ can be expressed by 
$\Omega^{k-2,0}=\int_{\gamma}\frac{\omega}{P}$, 
where $\omega=\sum_{i=1}^{k}(-1)^{i+1}dX_1\wedge ... \wedge \hat{dX_i} \wedge...\wedge dX_{k}$ and 
$\gamma$ is the small tube around the locus $\{P=0\}$. 
It is known that these periods are solutions of the Picard-Fuchs equation (\ref{hyper2}).

In the context of the B-model side of open mirror symmetry, the following integral over a $(k-2)$-chain $\Gamma^{k-2}$ is important: 
\begin{equation}
{\cal T}(z)=\int_{\Gamma^{k-2}}\Omega^{k-2,0}. \label{chainintegral}
\end{equation}
In this Appendix, we call (\ref{chainintegral}) the domainwall tension for the physical reason \cite{walcher,OV}. 
The boundaries of $(k-2)$-chain $\Gamma^{k-2}$ are the locus of B-brane. 
Mathematically, the formula of the domainwall tension (\ref{chainintegral}) is known as the normal function \cite{mw}. 
In our case, we take the following two boundaries ($(k-3)$-dimensional cycles) by generalizing the situation appearing in \cite{mw}. 
First, we take $(k-2)$-hyperplanes, $Q_1$, ..., $Q_{k-2}$, as follows: 
\begin{equation}
Q_1:X_1+X_2=0, \ \ \ 
Q_2:X_3+X_4=0, \ \ \ ..., \ \ \ 
Q_{k-2}:X_{k-2}+X_{k-1}=0. 
\end{equation}
Then we consider the complete intersection of $W_k^k$ and $Q_i$: 
\begin{equation}
P|_{Q_1,Q_2,...,Q_{k-2}}=X_k(X_k^{\frac{k-1}{2}}-\sqrt{k\psi}X_1X_3...X_{k-2})(X_k^{\frac{k-1}{2}}+\sqrt{k\psi}X_1X_3...X_{k-2}). 
\end{equation}
Boundaries are defined by two irreducible components $C_\pm:\{X_k^{\frac{k-1}{2}}\pm\sqrt{k\psi}X_1X_3...X_{k-2}=0\}$.\footnote{
$C_\pm$ and $RP^{k-2}$ are expected to be the mirror pair branes 
by applying (the higher-dimensional extended version of) the methods discussed in \cite{mw} or \cite{AV,AHMM}. }
For the mirror quintic case ($k=5$), it was observed that the above integral (\ref{chainintegral}) is a solution of the inhomogeneous version of the Picard-Fuchs equations \cite{walcher} and the inhomogeneous term was computed in \cite{mw}. 
An alternative method to evaluate the integral (\ref{chainintegral}) directly was developed \cite{FNSS} and applied to the pfaffian cases \cite{SS}. 
In the following, this direct integration method is generalized and applied to higher dimensional Calabi-Yau hypersurfaces. 
\\

First of all, 
we take local coordinates $(T,Y_1,...,Y_{\frac{k-3}{2}},\zeta_1,...,\zeta_{\frac{k-3}{2}},W)$ on the patch $\{X_1\neq 1\}$ as follows\footnote{
Contributions from the other patches result in just the overall factor of the final results. 
}: 
\begin{eqnarray}
&&X_1=1, X_2=T, \ \ \ 
X_3=Y_1^{\frac{k-1}{2k}}\zeta_1, X_4=Y_1^{\frac{k-1}{2k}}\zeta_1^{-1}, \ \ \ ..., \ \ \ 
X_{k-2}=Y_{\frac{k-3}{2}}^{\frac{k-1}{2k}}\zeta_{\frac{k-3}{2}}, 
X_{k-1}=Y_{\frac{k-3}{2}}^{\frac{k-1}{2k}}\zeta_{\frac{k-3}{2}}^{-1}, \\
&&X_k=(k\psi)^{\frac{1}{k-1}}W^{-\frac{2}{k-1}}(W^2-1)^{\frac{1}{k-1}}T^{\frac{1}{k-1}}Y_1^{\frac{1}{k}}...Y_{\frac{k-3}{2}}^{\frac{1}{k}}. 
\end{eqnarray}
The $W$ coordinate is introduced by $\frac{X_{k}^{k}}{(k\psi)X_1X_2...X_{k}}=\frac{W^2-1}{W^2}$. 
The form of $W$ is chosen for later use and the phase factor is ignored. 
Intuitively, $Y_1$, ..., $Y_{\frac{k-3}{2}}$ parametrize local coordinates on $C_\pm$, 
$T$ and $W$ are the coordinates which are longitudinal to $C_\pm$, and $\zeta_1$, ..., $\zeta_{\frac{k-3}{2}}$ are {\it polar} coordinates around $C_\pm$. 
These local coordinates are useful since after $Y_i$-integrations the period integral is factorized with respect to the remaining local coordinates. 
In the following, we will see this explicitly. 

The defining equation in the new coordinates becomes 
\begin{eqnarray}
P&=&1+T^k+Y_1^{\frac{k-1}{2}}f(\zeta_1) +...+Y_{\frac{k-3}{2}}^{\frac{k-1}{2}}f(\zeta_{\frac{k-3}{2}})
-z^{-\frac{1}{k-1}}(W^2-1)^{\frac{1}{k-1}}W^{-\frac{2k}{k-1}}T^{\frac{k}{k-1}}Y_1...Y_{\frac{k-3}{2}}\no\\
&=&1+T^{2m+1}+Y_1^{m}f(\zeta_1) +...+Y_{m-1}^{m}f(\zeta_{m-1})
-z^{-\frac{1}{2m}}(W^2-1)^{\frac{1}{2m}}W^{-\frac{1}{m}-2}T^{\frac{1}{2m}+1}Y_1...Y_{m-1}, 
\end{eqnarray}
where we introduce a symbol $f(\zeta):=\zeta^k+\zeta^{-k}$ for simplicity. 
In the second equality, we rewrite the formula by $m\in Z_{\geq 0}$ that obeys $k=2m+1$ and the genuine moduli $z:=\frac{1}{(k\psi)^k}$. 

The (normalized) period integrals in the new coordinates can be expressed by 
\begin{eqnarray}
k^{k-2}(k\psi)\int\frac{dX_2dX_3...dX_k}{(2\pi i)^{k-1}}\frac{1}{P}
=\frac{(2m)^{m-1}(2m+1)^{m}}{mz^{\frac{1}{2m}}}
\int\frac{dTdWdY_1dY_2...dY_{m-1}d\zeta_1 d\zeta_2 ... d\zeta_{m-1}}{(2\pi i)^{2m}\zeta_1\zeta_2 ... \zeta_{m-1}T^{-\frac{1}{2m}}W^{\frac{1}{m}+1}(W^2-1)^{-\frac{1}{2m}+1}}\frac{1}{P}.
\end{eqnarray}
The pre-factor is chosen by the usual convention. 
$\int$ symbolically means the integrals including both over $(k-2)$-cycles (or chain), and with respect to $\gamma$. 
All of them are performed as certain residue integrations. 
The difference between the fundamental period and the domainwall tension will be explained later. 

Now, we perform $Y_i$-integrals by picking up suitable poles. 
First we perform $Y_1$-integration. 
We regard $P$ as the degree $m$ polynomial with respect to $Y_1$ as follows: 
\begin{eqnarray}
P=A+z^{-\frac{1}{2m}}BY_1+CY_1^m, 
\end{eqnarray}
where 
\begin{eqnarray}
A&=&1+T^{2m+1}+Y_2^{m}f(\zeta_2)+...+Y_{m-1}^{m}f(\zeta_{m-1}), \\
B&=&-(W^2-1)^{\frac{1}{2m}}W^{-\frac{1}{m}-2}T^{\frac{1}{2m}+1}Y_2...Y_{m-1}, \\
C&=&f(\zeta_1). 
\end{eqnarray}
Near the large complex structure limit $z=0$, we can consider the following expansion: 
\begin{eqnarray}
\frac{1}{P}
&=&\frac{1}{A+z^{-\frac{1}{2m}}BY_1+CY_1^m}
=\frac{1}{(z^{-\frac{1}{2m}}B)(Y_1+\frac{A}{B}z^{\frac{1}{2m}})
\left[1+z^{\frac{1}{2m}}\frac{CY_1^m}{B(Y_1+\frac{A}{B}z^{\frac{1}{2m}})}\right]}\no\\
&=&z^{\frac{1}{2m}}\sum_{n=0}^{\infty}z^{\frac{n}{2m}} \frac{(-1)^nC^nY_1^{mn}}{B^{n+1}(Y_1+\frac{A}{B}z^{\frac{1}{2m}})^{n+1}}
=z^{\frac{1}{2m}}\oint\frac{ds}{2\pi i}\frac{\pi\cos(\pi s)}{\sin(\pi s)}z^{\frac{s}{2m}} \frac{(-1)^sC^sY_1^{ms}}{B^{s+1}(Y_1+\frac{A}{B}z^{\frac{1}{2m}})^{s+1}}. 
\end{eqnarray}
In the final equality, 
we perform an analytic continuation of the sum with respect to $n$ to the residue integral with respect to complex number $s$, 
where the $s$-integral is performed by a contour encircling the non-negative real axis. 
We expect that possible divergences can be avoided by this prescription. 
We can perform $Y_1$-integration by picking up a pole at $Y_1=-\frac{A}{B}z^{\frac{1}{2m}}$ 
and the period integral becomes 
\begin{eqnarray}
&&
\frac{(2m)^{m-1}(2m+1)^{m}}{m}
\oint\frac{ds}{2\pi i}\frac{\pi\cos(\pi s)}{\sin(\pi s)}z^{\frac{s}{2}}\frac{(-1)^{ms}\Gamma(ms+1)}{\Gamma(s+1)\Gamma((m-1)s+1)}\times\no\\
&&\times\int\frac{dTdWdY_2...dY_{m-1}d\zeta_1 d\zeta_2 ... d\zeta_{m-1}}{(2\pi i)^{2m-1}\zeta_1\zeta_2 ... \zeta_{m-1}T^{-\frac{1}{2m}}W^{\frac{1}{m}+1}(W^2-1)^{-\frac{1}{2m}+1}}\frac{A^{(m-1)s}C^s}{B^{ms+1}}\no\\
&=&
\frac{(2m)^{m-1}(2m+1)^{m}}{m}
\oint \frac{ds}{2\pi i} \frac{\pi\cos(\pi s)}{\sin(\pi s)}z^{\frac{s}{2}}\frac{(-1)\Gamma(ms+1)}{\Gamma(s+1)\Gamma((m-1)s+1)}\int\frac{dTdWdY_2...dY_{m-1}d\zeta_1 d\zeta_2 ... d\zeta_{m-1}}{(2\pi i)^{2m-1}\zeta_1\zeta_2 ... \zeta_{m-1}}\times\no\\
&&\times \frac{f(\zeta_1)^s}{(W^2-1)^{\frac{s}{2}+1}W^{-2ms-s-1}T^{ms+\frac{s}{2}+1}}
\frac{[{1+T^{2m+1}+Y_2^{m}f(\zeta_2)+...+Y_{m-1}^{m}f(\zeta_{m-1})}]^{(m-1)s}}{{Y_2^{ms+1}...Y_{m-1}^{ms+1}}}. 
\end{eqnarray}
Then we perform $Y_2$, $Y_3$, ... $Y_{m-1}$-integrations in turn. 
By picking up poles at $Y_i=0$ $(i=2$, ..., $m-1)$, 
we obtain the following factorized integral formula 
\footnote{The fact that the period integral is expanded by $z^{\frac{1}{2}}$ at this stage is related to $Z_2$ open moduli structure. }
: 
\begin{eqnarray}
&&
\frac{(2m)^{m-1}(2m+1)^{m}}{m}
\oint \frac{ds}{2\pi i} \frac{\pi\cos(\pi s)}{\sin(\pi s)}z^{\frac{s}{2}}\frac{(-1)^{\frac{s}{2}}\Gamma(ms+1)}{\Gamma(s+1)\Gamma((m-1)s+1)}\int\frac{dTdWd\zeta_1 d\zeta_2 ... d\zeta_{m-1}}{(2\pi i)^{m+1}\zeta_1\zeta_2 ... \zeta_{m-1}}\times\no\\
&&\times \frac{f(\zeta_1)^s}{(1-W^2)^{\frac{s}{2}+1}W^{-2ms-s-1}T^{ms+\frac{s}{2}+1}}
\frac{\Gamma((m-1)s+1)}{\Gamma(s+1)^{m-1}}(1+T^{2m+1})^sf(\zeta_2)^s...f(\zeta_{m-1})^s\no\\
&=&
\frac{(2m)^{m-1}(2m+1)^{m}}{m}
\oint\frac{ds}{2\pi i}\frac{\pi\cos(\pi s)}{\sin(\pi s)}z^{\frac{s}{2}}\frac{(-1)^{\frac{s}{2}}\Gamma(ms+1)}{\Gamma(s+1)^m}\times\no\\
&&\times\left(\int\frac{d\zeta}{(2\pi i)\zeta}f(\zeta)^s\right)^{m-1}
\int \frac{dT}{2\pi i}\frac{(1+T^{2m+1})^s}{T^{ms+\frac{s}{2}+1}}
\int \frac{dW}{2\pi i}\frac{(1-W^2)^{-\frac{s}{2}-1}}{W^{-2ms-s-1}}. 
\end{eqnarray}

The remaining task is to compute the remaining integrals individually. 
We proposed the following simple prescription in \cite{FNSS,SS}: 
For the fundamental period $w_0$ (the unique period which is regular under expansion around $z=0$), 
we choose contour integrals for all integrals, and 
for the domainwall tension ${\cal T}$, we choose a line integral for $W$-integral and contour integrals for the other integrals. 
More precisely, in the present situation, we take the following integral paths according to the orbifold structure. 
For the fundamental period, 
the path of the $\zeta$-integral is $\zeta=e^{i\phi}$ $(-\frac{\pi}{2}\leq \phi\leq -\frac{\pi}{2}+\frac{\pi}{2m+1})$, 
the path of the $T$-integration is $T=e^{i\theta}$ $(-\pi\leq \theta\leq -\pi+\frac{2 \pi}{2m+1})$, and 
the path of the $W$-integration is the contour encircling $-1$ and $1$ one times. 
For the domainwall tension, 
we take a line integral from $-1$ to $1$ for the $W$-integral (note that B-brane is expressed by $W=\pm1$) and 
the other integrals is performed by the same paths as those of the fundamental period. 
Moreover, the tension of the domainwall between two branes $C_+$ and $C_-$ is expressed by two parts (contributions from $C_+$ and $C_-$) as follows: 
\begin{equation}
{\cal T}=\int_{-1}^{1}dW(...)={\cal T}_+-{\cal T}_-=\int_{0}^{1}dW(...)+\int_{-1}^{0}dW(...). 
\end{equation}
Since it turns out that ${\cal T}_-(z^{\frac{1}{2}})={\cal T}_+(-z^{\frac{1}{2}})$, 
we can concentrate on the evaluation of ${\cal T}_+$ in the following. 

Each integral can be evaluated as follows: 
\begin{eqnarray}
&&\int\frac{d\zeta}{(2\pi i)\zeta}f(\zeta)^s =\frac{(-1)^{ms}}{2(2m+1)}\frac{\Gamma(s+1)}{\Gamma(\frac{s}{2}+1)^2},\\
&&\int \frac{dT}{2\pi i}\frac{(1+T^{2m+1})^s}{T^{ms+\frac{s}{2}+1}}=\frac{(-1)^{ms}}{2m+1}
\frac{\Gamma(s+1)}{\Gamma(\frac{s}{2}+1)^2}, \\
&&\oint \frac{dW}{2\pi i}\frac{(1-W^2)^{-\frac{s}{2}-1}}{W^{-2ms-s-1}}=(-1)^{ms}\cos((2m+1)\frac{\pi s}{2})\frac{\Gamma((2m+1)\frac{s}{2}+1)}{\Gamma(ms+1)\Gamma(\frac{s}{2}+1)}, \\
&&\int_0^{1} \frac{dW}{2\pi i}\frac{(1-W^2)^{-\frac{s}{2}-1}}{W^{-2ms-s-1}}
=-\frac{1}{4i}\frac{1}{\sin(\frac{\pi s}{2})}\frac{\Gamma((2m+1)\frac{s}{2}+1)}{\Gamma(ms+1)\Gamma(\frac{s}{2}+1)}. 
\end{eqnarray}

The fundamental period $w_0$ becomes 
\begin{eqnarray}
w_0(z)=m^{m-2}\oint\frac{ds}{2\pi i}\frac{\pi\cos(\pi s)\cos((2m+1)\frac{\pi s}{2})}{\sin(\pi s)}
z^{\frac{s}{2}}
\frac{\Gamma((2m+1)\frac{s}{2}+1)}{\Gamma(\frac{s}{2}+1)^{2m+1}}
(-1)^{\frac{s}{2}}, 
\end{eqnarray}
and the domainwall tension becomes 
\begin{eqnarray}
{\cal T}_+(z)=\frac{m^{m-2}}{4i}\oint\frac{ds}{2\pi i}\frac{\pi\cos(\pi s)}{\sin(\pi s)\sin(\frac{\pi s}{2})}
z^{\frac{s}{2}}
\frac{\Gamma((2m+1)\frac{s}{2}+1)}{\Gamma(\frac{s}{2}+1)^{2m+1}}(-1)^{m^2s+\frac{s}{2}}. 
\end{eqnarray}

Finally, we perform the $s$-integral by picking up appropriate poles and 
rewrite by $2m+1=k$ and $z=e^x$. 
We neglect the common overall factor $m^{m-2}=\left(\frac{k-1}{2}\right)^{\frac{k-5}{2}}$ in the following. 
For the fundamental period, we have simple poles at $s=2d$ $(d\in Z_{\geq 0})$ and obtain the well-known formula: 
\begin{equation}
w_0(x)=\sum_{d=0}^{\infty}\frac{\Gamma(kd+1)}{\Gamma(d+1)^k}e^{dx}. 
\label{F}
\end{equation}
For ${\cal T}_\pm$, 
we have simple poles at $s=2d+1$ $(d\in Z_{\geq 0})$ and double poles at $s=2d$ $(d\in Z_{\geq 0})$, and obtain \footnote{
Here, $A=1$ for even $m$ and $A=3$ for odd $m$. 
At this stage, we haven't clarified the meaning of this constant. 
}
\begin{equation}
{\cal T}_\pm(x)=\frac{1}{4\pi i}w_1(x)\pm \frac{A}{4}w_0(x)\pm \frac{1}{8}\left(\frac{2}{\pi}\right)^{\frac{k-1}{2}}\tau_k(x), 
\label{T}
\end{equation}
where $w_0$ is the fundamental period (\ref{F}), $w_1$ is the logarithmic period, and $\tau_k$ is given in (\ref{superpotential}): 
\begin{equation}
\tau_{k}(x):=\sum_{d=1}^{\infty}2\frac{(k(2d-1))!!}{((2d-1)!!)^{k}}e^{\frac{2d-1}{2}x}. 
\end{equation}
Thus we obtain (\ref{superpotential}) as a part of the domainwall tension, 
not by solving (\ref{hyper3}) but by performing the integration directly. 
Contributions from the fundamental and logarithmic period in (\ref{T}) are important when discussing the monodromy property \cite{walcher}. 

We expect that this method can be generalized to the cases of projective complete intersections and Fano projective hypersurfaces. 
It is known that off-shell open mirror symmetry has quite rich structures \cite{JS,AHMM}. 
Application of the direct integration method to off-shell cases also seems to be possible as was done for some Calabi-Yau $3$-fold cases \cite{FNSS}.


\end{document}